\documentclass{compositio}

\usepackage{amsmath}

\newtheorem{Theorem}{Theorem}[section]
\newtheorem{Proposition}[Theorem]{Proposition}
\newtheorem{Corollary}[Theorem]{Corollary}
\newtheorem{Lemma}[Theorem]{Lemma}
\newtheorem{Definition}[Theorem]{Definition}
\theoremstyle{remark}
\newtheorem{Remark}[Theorem]{Remark}

\numberwithin{equation}{section}

\def\ta{\theta}
\def\sgn{\operatorname{sgn}}

\begin{document}

\title[Elliptic determinant evaluations and Macdonald identities]
{Elliptic  determinant evaluations
and the Macdonald identities for affine root systems}

\author{Hjalmar Rosengren}
\email{hjalmar@math.chalmers.se}
\address{Department of Mathematical Sciences, Chalmers University of Technology and
G\"ote\-borg University, SE-412 96 G\"oteborg, Sweden}


\author{Michael Schlosser}
\email{schlosse@ap.univie.ac.at}
\address{Fakult\"at f\"ur Mathematik der Universit\"at Wien,
Nordbergstra{\ss}e 15, A-1090 Wien, Austria}


\thanks{Research partially supported by EC's IHRP Programme,
grant HPRN-CT-2001-00272 ``Algebraic Combinatorics in Europe''. The
first author was  supported by the Swedish Science Research
Council (VR).
The second author was fully supported by
FWF Austrian Science Fund grant \hbox{P17563-N13}}


\begin{abstract}
We obtain several determinant evaluations, related to affine root
systems, which provide
 elliptic extensions of  Weyl denominator formulas.
 Some of these are new, also in the polynomial
special case, while others yield new proofs of the
 Macdonald identities  for the seven
infinite families of irreducible reduced affine root systems. 
\end{abstract}

\maketitle

\section{Introduction}

Determinants play an important role in many areas of mathematics.
Often, the solution of a particular problem in combinatorics,
mathematical physics  or, simply, linear algebra, depends on the explicit
computation of a determinant.
Some useful and efficient tools for evaluating
determinants are provided in Krattenthaler's survey articles
\cite{K2}, \cite{K3},
which also contain many explicit determinant evaluations that have
appeared in the literature
 and give references where further such formulae can be found.

As  examples of interesting  determinant evaluations, we mention  the Weyl
denominator formulas for classical root systems, which play a
fundamental role  in Lie theory and related areas. In general, the
Weyl denominator formula for a reduced root system reads
\begin{equation}\label{gwd}\sum_{w\in W}\det(w)e^{w(\rho)-\rho}=\prod_{\alpha\in R_+}(1-e^{-\alpha}), \end{equation}
where $W$ is the Weyl group, $R_+$ the set of positive roots and
$\rho=\frac 12\sum_{\alpha\in R_+}\alpha$. For the classical
  root systems $A_{n-1}$, $B_n$, $C_n$
and $D_n$, this identity takes the explicit form
\begin{subequations}\label{wd}
\begin{align}
\label{awd}
\det_{1\leq i,j\leq n}\left(x_i^{j-1}\right)
&=\prod_{1\leq i<j\leq n}(x_j-x_i),\\
\label{bwd}
\det_{1\leq i,j\leq n}\left(x_i^{j-n}-x_i^{n+1-j}\right)
&=\prod_{i=1}^nx_i^{1-n}(1-x_i)\prod_{1\leq i<j\leq
  n}(x_j-x_i)(1-x_ix_j), \\
\label{cwd}
\det_{1\leq i,j\leq n}\left(x_i^{j-n-1}-x_i^{n+1-j}\right)
&=\prod_{i=1}^nx_i^{-n}(1-x_i^2)\prod_{1\leq i<j\leq
  n}(x_j-x_i)(1-x_ix_j),\\
\label{dwd}
\det_{1\leq i,j\leq n}\left(x_i^{j-n}+x_i^{n-j}\right)&
=2\prod_{i=1}^nx_i^{1-n}\prod_{1\leq i<j\leq n}(x_j-x_i)(1-x_ix_j), 
\end{align}
\end{subequations}
respectively.

In this article, we are interested in generalizing \eqref{wd} to the level of 
\emph{elliptic determinant evaluations}. By this we  mean that  the
matrix elements should be defined in terms of theta functions, so that it is \emph{a priori}
clear that the quotient of the two sides of the identity is an
elliptic function of some natural parameters. 
 Up to date, according to our
knowledge,  very few elliptic determinant  (and  pfaffian)
evaluations are known, see \cite{Fr}, \cite{FS}, \cite[Lem.~1]{H}, \cite{O}, 
\cite[Th.\ 2.10]{R2},
\cite[App.~B]{TV} and
\cite[Th.~4.17, Lem.~5.3]{W}. Most of these results contain elliptic
extensions of Weyl denominators, and are thus apparently related to root systems.

An  elliptic extension of the Weyl denominator
formula was obtained by Macdonald \cite{M}, see also \cite{D}. He
introduced, and completely classified, \emph{affine
  root systems}.
Moreover, he  extended  the  Weyl denominator formula to the case of
 \emph{reduced} affine root systems. In this setting, both the root
 system and the Weyl group are infinite, so  the resulting
 \emph{Macdonald identities} equate an infinite series and an infinite
 product. The precise statement is  more complicated than
 \eqref{gwd}, see \cite[Th.~8.1]{M} and, for the special cases of
 interest to us, Corollary \ref{mlc} below.
The Macdonald
 identities can be interpreted in terms of 
 Kac--Moody algebras \cite{Ka}.  Notable special cases 
include Watson's quintuple product identity \cite{Wa} (for the affine
root system $BC_1$), Winquist's identity \cite{Wi} (for  $B_2$)
and the so called  septuple product identity \cite{fk,hi1,hi2}
(put $x_2=-1$ in the $BC_2$ case of
Proposition \ref{mdp} below).

There are seven infinite families of
irreducible reduced affine root systems  and seven exceptional
cases. We will only consider  the  infinite families, which
Macdonald denotes $A$,  
$B$, $B^\vee$, $C$, $C^\vee$, $BC$ and $D$. 
They should not be confused with the  classical root systems mentioned
above. (For instance, the classical root system $BC_n$ is non-reduced
whereas the affine root system $BC_n$ is reduced.)
Although the corresponding  Macdonald
   identities do give  elliptic extensions of
\eqref{wd}, it is only for type
  $C$, $C^\vee$ and $BC$ that they can immediately
be written as determinant evaluations. Nevertheless, one of our goals 
is to rewrite all seven cases in determinant form, and  prove them by an
``identification of factors'' argument similar to the usual proof
of the Vandermonde determinant \eqref{awd}.
This new proof of the Macdonald
 identities is rather similar to Stanton's elementary proof \cite{St},
 but the use of determinants makes the details more streamlined.

For each  affine root system $R$ under consideration, we define a
corresponding notion of $R$ theta function. We then give a ``master
determinant formula'',  Proposition~\ref{wp}, which expresses a
determinant of $R$ theta functions as a constant times the $R$
 Macdonald denominator. When the constant can be
 explicitly determined, we have a genuine determinant evaluation. 
Such explicit instances of the master formula include a
 determinant of Warnaar (Proposition \ref{bcdet} below), new
 generalized Weyl
 denominator formulas for all seven families of reduced affine root systems (Theorem~\ref{adet}, Theorem~\ref{cdet} and Corollaries \ref{bcdet2}  to
 \ref{ddet}) and determinant versions of the Macdonald identities 
(Proposition~\ref{mdp}). Theorem \ref{adet} include as
 special cases the determinants of  Frobenius and Hasegawa cited
 above, and has a non-trivial overlap with the determinant of Tarasov
 and Varchenko.

The most striking difference between our new elliptic denominator formulas and
 those found by Macdonald
 is the  large number of free
parameters in our identities. This probably makes the results more difficult to interpret 
 in terms of, say,
affine Lie algebras. On the other hand,  the presence of free parameters
seems useful for certain applications. Indeed, special cases of our identities have found
applications to {\em multidimensional basic and elliptic
hypergeometric series and integrals}, see  \cite{GK}, \cite{KN}, \cite{R},
\cite{R2}, \cite{Ro1},
\cite{Ro},  \cite{RS1},  \cite{S1},
\cite{S2}, \cite{S3},  \cite{S4}, \cite{Sp},   \cite{W}, to the study of
\emph{Ruijsenaars operators} and related integrable systems 
\cite{H}, \cite{Ru}, to \emph{combinatorics}, see \cite{K2} for an
extensive list of references,
as well as   to {\em number theory} \cite{Ro2}.
It  thus seems very
likely that our new results will find similar applications.

Our paper is organized as follows. 
Section~\ref{secpre} contains preliminaries on Jacobi theta functions. In
Section~\ref{sectrs} we introduce theta functions associated to the
seven families of reduced affine root
systems.  We then
give our master formula, Proposition~\ref{wp}. In
Section~\ref{secell} we obtain several elliptic determinant evaluations
that can be viewed as explicit versions of  Proposition \ref{wp}.
The main results  are  Theorems~\ref{adet} and~\ref{cdet} (the other determinant evaluations are corollaries of
these).
 Section~\ref{secrat} features several
corollaries obtained by restricting to the {\em polynomial} special
case. Finally, in Section~\ref{secmac}, we obtain  determinant
evaluations
that are shown to be equivalent to the Macdonald identities for 
non-exceptional reduced affine root systems.

\begin{acknowledgements}
 We thank Eric Rains for his comments on Corollary
\ref{lem}, leading to some improvements in that part of the paper,
 and Vitaly Tarasov for clarifying how to
obtain Corollary \ref{tvcor} from the results of
 \cite{TV}, see Remark \ref{tvrem}.
\end{acknowledgements}

\section{Preliminaries}\label{secpre}

Throughout this paper, we implicitly assume that all scalars are generic,
so that no denominators in our identities vanish.

The letter $p$ will denote
 a fixed number such that  $0<|p|<1$. When dealing with the root
 system $C_n^\vee$, we will also assume a fixed choice of square root
 $p^{1/2}$. The  case  $p=0$ will be considered  in
 Section \ref{secrat}.

We  use the standard notation
$$(a)_\infty=(a;p)_\infty=\prod_{j=0}^\infty(1-ap^j),$$
$$(a_1,\dots,a_n)_\infty=(a_1,\dots,a_n;p)_\infty=(a_1;p)_\infty\dotsm(a_n;p)_\infty. $$
Then, 
\begin{equation}\label{qr}
(x^k;p^k)_\infty=\prod_{j=0}^{k-1}(x\omega_k^j;p)_\infty,
\qquad
(x;p)_\infty=\prod_{j=0}^{k-1}(xp^j;p^k)_\infty,
\end{equation}
where $\omega_k$ denotes a primitive $k$th root of unity.

We employ ``multiplicative'', rather than
``additive'', notation for theta functions. 
This corresponds to realizing the torus
$\mathbb C/(\mathbb Z+\tau\mathbb Z)$ as $(\mathbb
C\setminus\{0\})/(z\mapsto pz)$, where $p=e^{2\pi i\tau}$. 
Thus, we  take as our building block the  function
$$\ta(x)=\ta(x;p)=(x,p/x;p)_\infty. $$
We will sometimes use the shorthand notation
$$\theta(a_1,\dots,a_n)=\theta(a_1)\dotsm\theta(a_n), $$
$$\theta(xy^\pm)=\theta(xy)\theta(x/y). $$

The function $\theta(x)$ is holomorphic for $x\neq 0$ and has single zeroes
precisely at $p^{\mathbb Z}$. 
 Up to an elementary factor,
 $\theta(e^{2\pi ix};e^{2\pi i \tau})$ equals
the Jacobi theta function $\theta_1(x|\tau)$. 
We will frequently use the inversion formula
$$\theta(1/x)=-\frac1 x\,\theta(x) $$
and the quasi-periodicity
$$\ta(px)=-\frac1 x\,\theta(x).$$
By Jacobi's triple product identity, we have the Laurent expansion
\begin{equation}\label{jti}
\theta(x)=\frac{1}{(p)_\infty}\sum_{n=-\infty}^\infty(-1)^np^{\binom
n2}x^n.\end{equation}

Similarly to \eqref{qr}, we have
\begin{equation}\label{qrt}
\theta(x^k;p^k)=\prod_{j=0}^{k-1}\theta(x\omega_k^j;p),
\qquad
\theta(x;p)=\prod_{j=0}^{k-1}\theta(xp^j;p^k),
\end{equation}
which, when $k=2$, implies
\begin{equation}\label{txsq}
\theta(x^2)=\theta(x,-x,p^\frac12 x,-p^\frac12 x).\end{equation}
Since $\theta(x)$ has a single zero at $x=1$, it follows that
\begin{equation}\label{tev}
\theta(-1,p^{\frac 12},-p^{\frac 12})=\lim_{x\rightarrow
  1}\frac{\theta(x^2)}{\theta(x)}=2. 
\end{equation}

\section{Theta functions on root systems}\label{sectrs}

The Macdonald identities involve the Macdonald denominator 
\begin{equation}\label{mcd}\prod_{\alpha\in R_+}(1-e^{-\alpha}),\end{equation}
 where
$R_+$ is the positive part of a reduced affine root system and
$e^\alpha$ a formal exponential. Although we
 will not need anything of Macdonald's theory,  it may be
instructive to explain what \eqref{mcd} means in the case $R=C_n$. 
Let $e_i$, $1\leq i\leq n$, be a basis for $\mathbb R^n$, and
write $k+\varepsilon_i$ for the affine function $e_j\mapsto
k+\delta_{ij}$. Then, affine  $C_n$ consists of the roots
$$k\pm 2\varepsilon_i,\qquad k\in\mathbb Z,\ 1\leq i\leq n,$$
$$k\pm \varepsilon_i\pm\varepsilon_j,\qquad k\in\mathbb Z,\ 1\leq
i<j\leq n.$$
The positive roots are
$$k+2\varepsilon_i,\qquad k\geq 0,\ 1\leq i\leq n, $$
$$k-2\varepsilon_i,\qquad k\geq 1,\ 1\leq i\leq n, $$
$$k+\varepsilon_i+\varepsilon_j,\quad k+\varepsilon_i-\varepsilon_j, \qquad k\geq 0,\ 1\leq
i<j\leq n,  $$
$$k-\varepsilon_i+\varepsilon_j,\quad k-\varepsilon_i-\varepsilon_j, \qquad k\geq 1,\ 1\leq
i<j\leq n.  $$
Thus, the Macdonald denominator for $C_n$ is
\begin{multline}\label{cmd}\prod_{i=1}^n\prod_{k=0}^\infty(1-e^{-(k+2\varepsilon_i)})(1-e^{-(k+1-2\varepsilon_i)})\\
\times\prod_{1\leq
i<j\leq n}\prod_{k=0}^\infty(1-e^{-(k+\varepsilon_i+\varepsilon_j)})
(1-e^{-(k+\varepsilon_i-\varepsilon_j)})(1-e^{-(k+1-\varepsilon_i+\varepsilon_j)})(1-e^{-(k+1-\varepsilon_i-\varepsilon_j)}).
 \end{multline} 
Introducing variables $p$
and $x_1,\dots,x_n$ by $p=e^{-1}$, $x_i=p^{-\varepsilon_i}$, \eqref{cmd}
takes the form
$$\prod_{i=1}^n\theta(x_i^2)\prod_{1\leq
i<j\leq n}\theta(x_ix_j^\pm),$$
where $\theta(x)=\theta(x;p)$.
The $C_n$ Macdonald identity gives the explicit multiple Laurent
expansion of this function, where $x_i$ are viewed as non-zero complex
variables and $p$ as a constant with $|p|<1$. 

More generally, the Macdonald denominators for the seven families of
reduced affine root systems equal, up to a trivial factor that has
been chosen for convenience,
\begin{align*}
W_{A_{n-1}}(x)&=\prod_{1\leq i<j\leq n}x_j\theta(x_i/x_j), \\
W_{B_n}(x)&=\prod_{i=1}^n\theta(x_i)
\prod_{1\leq i<j\leq n}x_i^{-1}\theta(x_ix_j^\pm), \\
W_{B_n^\vee}(x)&=\prod_{i=1}^n x_i^{-1}\theta(x_i^2;p^2)
\prod_{1\leq i<j\leq n}x_i^{-1}\theta(x_ix_j^\pm), \\
W_{C_n}(x)&=\prod_{i=1}^n x_i^{-1}\theta(x_i^2)
\prod_{1\leq i<j\leq n}x_i^{-1}\theta(x_ix_j^\pm),\\
W_{C_n^\vee}(x)&=\prod_{i=1}^n\theta(x_i;p^{\frac 12})
\prod_{1\leq i<j\leq n}x_i^{-1}\theta(x_ix_j^\pm), \\
W_{BC_n}(x)&=\prod_{i=1}^n\theta(x_i)\theta(px_i^2;p^2)
\prod_{1\leq i<j\leq n}x_i^{-1}\theta(x_ix_j^\pm), \\
W_{D_n}(x)&=\prod_{1\leq i<j\leq n}x_i^{-1}\theta(x_ix_j^\pm).
\end{align*}
We will use the above list
as a rule for labelling our results. Each of our elliptic determinant
evaluations expresses the Macdonald denominator of some affine root system
as a determinant.

The following definition may seem strange, since root systems are usually
associated to multivariable functions. However, it will enable
us to give a very succinct statement of Proposition \ref{wp}. Note
that, except in the case $R=A_{n-1}$,
$W_R$ is an $R$ theta function
of each $x_i$.  This is easy to check directly, and is also clear
from Proposition  \ref{wp}.

\begin{Definition}\label{rtd}
 Let $f(x)$ be  holomorphic for $x\neq
  0$. Then, we call $f$ an $A_{n-1}$ theta function of norm $t$ if
\begin{equation}\label{ade}f(px)=\frac{(-1)^n}{tx^n}\,f(x).\end{equation}
Moreover, if $R$ denotes either  $B_n$,
$B^\vee_n$, $C_n$, $C^\vee_n$, $BC_n$ or $D_n$, 
  we call $f$ an $R$ theta function if
\begin{align*}
f(px)&=-\frac{1}{p^{n-1}x^{2n-1}}\,f(x),& f(1/x)&=-\frac
1x\,f(x),& R&=B_n,\\
f(px)&=-\frac{1}{p^{n}x^{2n}}\,f(x),& f(1/x)&=-f(x),&
R&=B_n^\vee,\\
f(px)&=\frac{1}{p^{n+1}x^{2n+2}}\,f(x),& f(1/x)&=-f(x),&
R&=C_n,\\
f(px)&=\frac{1}{p^{n-\frac12}x^{2n}}\,f(x),& f(1/x)&=-\frac
1x\,f(x),& R&=C_n^\vee,\\
f(px)&=\frac{1}{p^{n}x^{2n+1}}\,f(x),& f(1/x)&=-\frac
1x\,f(x),& R&=BC_n,\\
f(px)&=\frac{1}{p^{n-1}x^{2n-2}}\, f(x),& f(1/x)&=f(x),&
R&=D_n.
\end{align*}
\end{Definition}

These notions depend on our fixed parameter $p$, and in the
case of $C_n^\vee$ on a choice of square root  $p^{1/2}$.

The following result gives  useful factorizations of $R$
theta functions. 

\begin{Lemma}\label{fl}
The function $f$ is an $A_{n-1}$ theta function of norm $t$ 
if and only if there exist constants $C$, $b_1,\dots,b_{n}$ such that
$b_1\dotsm b_n=t$ and
$$f(x)=C\,\theta(b_1x,\dots,b_nx). $$
For the other six cases, $f$ is an 
$R$ theta function if and only if there exist
constants $C$, $b_1,\dots,b_{n-1}$ such that
\begin{align*}f(x)&=C\,\ta(x)\,\ta(b_1x^\pm,\dots,b_{n-1}x^\pm),&
  R&=B_n,\\
f(x)&=C\,x^{-1}\ta(x^2;p^2)\,\ta(b_1x^\pm,\dots,b_{n-1}x^\pm),&
R&=B_n^\vee,\\
f(x)&=C\,x^{-1}\ta(x^2)\,\ta(b_1x^\pm,\dots,b_{n-1}x^\pm),& R&=C_n,\\
f(x)&=C\,\ta(x;p^{\frac12})\,\ta(b_1x^\pm,\dots,b_{n-1}x^\pm),&
R&=C_n^\vee,\\
f(x)&=C\,\ta(x)\ta(px^2;p^2)\,\ta(b_1x^\pm,\dots,b_{n-1}x^\pm),&
R&=BC_n,\\
f(x)&=C\,\ta(b_1x^\pm,\dots,b_{n-1}x^\pm),& R&=D_n,
\end{align*}
where $\theta(x)=\theta(x;p)$.
\end{Lemma}

\begin{proof}
Up to the change of variable $x\mapsto e^{2\pi ix}$, what we call
an $A_{n-1}$ theta function is usually called a theta function of
order $n$. In that case, the factorization theorem is classical, see
\cite[p.\ 45]{We}. Nevertheless, we review the proof. The ``if'' part is straight-forward, so we assume that $f$ is an $A_{n-1}$
theta function. Let $N$ be the number of zeroes of $f$, counted with
multiplicity,  inside any
period annulus $A=\{|p|r< |x|\leq r\}$. It is well-known that
$$N=\int_{\partial A}\frac{f'(x)}{f(x)}\,\frac{dx}{2\pi i}.$$
The equality \eqref{ade} differentiates to
$$\frac{f'(x)}{f(x)}-p\frac{f'(px)}{f(px)}=\frac{n}{x}, $$
which gives $N=n$. Thus, there exist $b_1,\dots, b_n$ so that the
zeroes, counted with multiplicity, are enumerated by $p^mb_i$,
$m\in\mathbb Z$, $i=1,\dots,n$. The function
$g(x)=f(x)/\theta(b_1x,\dots,b_nx)$ is then analytic for $x\neq 0$
and satisfies $g(px)=g(x)$, so by Liouville's theorem it is
constant. 
Finally, if  $f$ has norm $t$, one  checks  that
$b_1\dotsm b_n=t$.

Let us now  consider the case $R=D_n$. 
Since any  $D_n$
theta function $f$ is an $A_{2n-3}$ theta function, it has $2n-2$ zeroes in
each period annulus. It is easy to check from the definition that if 
 $a$ is a zero, then $1/a$ is a  zero of the same multiplicity, and if
 some zero should satisfy $a^2\in p^\mathbb Z$, then its multiplicity is
even. Thus, there exist $a_1,\dots,a_{n-1}$ so that the zeroes, with
multiplicity, are enumerated by $ p^m a_i^\pm$, $m\in\mathbb Z$,
$i=1,\dots, n-1$. As before, 
$g(x)=f(x)/\theta(a_1x^\pm,\dots,a_{n-1}x^\pm)$ is  analytic for $x\neq 0$
and satisfies $g(px)=g(x)$, so by Liouville's theorem it is constant.

The other cases are easily deduced from the case $R=D_n$. For
instance, assume that $f$ is a $BC_n$ theta function. Letting
$x=1$, $x=1/\sqrt p$ and $x=-1/\sqrt p$ in Definition \ref{rtd}, one finds that 
$f$ vanishes at these points and thus  
 $f(p^m)=f(\pm\sqrt p p^m)=0$ for any $m\in\mathbb Z$.  It follows that
 $g(x)=f(x)/\theta(x)\theta(px^2;p^2)$ is analytic for $x\neq 0$. It is
 straight-forward to check that $g$ is a $D_n$ theta function, so the
 desired factorization follows from the case $R=D_n$. The remaining cases can
 be treated similarly.
\end{proof}
 
We will also use the following result, which
 expresses  $R$ theta functions, when $R$ is not of type $A$, 
in terms of type $A$ theta functions.

\begin{Lemma}\label{dl}
The function $f$ is an $R$ theta function if and only if there exists
a function $g(x)$, holomorphic for $x\neq 0$, such that 
\begin{align*}
g(px)&=-\frac{1}{p^{n-1}x^{2n-1}}\,g(x),& f(x)&=g(x)-xg(1/x),
& R&=B_n,\\
g(px)&=-\frac{1}{p^{n}x^{2n}}\,g(x),& f(x)&=g(x)-g(1/x),&
R&=B_n^\vee,\\
g(px)&=\frac{1}{p^{n+1}x^{2n+2}}\,g(x),& f(x)&=g(x)-g(1/x),&
R&=C_n,\\
g(px)&=\frac{1}{p^{n-\frac12}x^{2n}}\,g(x),& f(x)&=g(x)-xg(1/x),
& R&=C_n^\vee,\\
g(px)&=\frac{1}{p^{n}x^{2n+1}}\,g(x),& f(x)&=g(x)-xg(1/x),& R&=BC_n,\\
g(px)&=\frac{1}{p^{n-1}x^{2n-2}}\, g(x),& f(x)&=g(x)+g(1/x),&
R&=D_n.
\end{align*}
\end{Lemma}

\begin{proof}
If $f$ is an $R$ theta function,  one may in each case choose
$g=f/2$. The converse is straight-forward.
\end{proof}

An important example, to be used later, is the case when
$R=C_1$ and $g(x)=x^{-2}\theta(ax,bx,cx,dx)$, $abcd=1$. Combining Lemma
\ref{fl} and Lemma \ref{dl} gives $g(x)-g(1/x)=C\,x^{-1}\theta(x^2)$, where $C$ may
be computed by plugging in $x=a$. This leads to the identity
\begin{equation}\label{radd}\frac
  1{x^{2}}\,\theta(ax,bx,cx,dx)-x^2\theta(a/x,b/x,c/x,d/x)
=\frac
1{ax}\,\theta(ab,ac,ad,x^2),\qquad abcd=1,\end{equation}
which is equivalent to Riemann's addition formula (cf.\ \cite[p.~451,
Example~5]{WW}).

We are now in a position to state our ``master formula''.

\begin{Proposition}\label{wp}
Let $f_1,\dots,f_n$ be  $A_{n-1}$ theta functions of norm
$t$. Then, 
\begin{subequations}\label{mf}
\begin{equation}\label{awpi}\det_{1\leq i,j\leq n}\left(f_j(x_i)\right)=C\,\theta(tx_1\dotsm x_n)
\,W_{A_{n-1}}(x)\end{equation}
for some constant $C$. Moreover, if $R$ denotes either  $B_n$,
$B^\vee_n$, $C_n$, $C^\vee_n$, $BC_n$ or $D_n$ and   $f_1,\dots,f_n$
are $R$ theta functions, we have
\begin{equation}\label{wpi}
\det_{1\leq i,j\leq n}\left(f_j(x_i)\right)=C\,W_R(x) 
\end{equation}
\end{subequations}
for some constant $C$.
\end{Proposition}

\begin{proof}
Consider first the case of \eqref{awpi}. For fixed $i=1,\dots,n$, 
let  $L(x_i)$ and $R(x_i)$
denote  the left-hand and right-hand sides, viewed as
functions of $x_i$.  It is straight-forward to verify that both $L$ and
$R$ are $A_{n-1}$ theta functions of  norm $t$. Thus, $f=L/R$ satisfies $f(px)=f(x)$, so if we can prove that
$f$ is analytic, it follows from Liouville's theorem that
it is constant. Up to multiplication with  $p^{\mathbb Z}$, the
zeroes of $R$ are situated at $x_i=x_j$, $j\neq i$ and at
$x_i=1/tx_1\dotsm 
\hat x_i\dotsm x_n$. For generic values of  $x_j$,
$j\neq i$, they are all single zeroes, so it is 
 enough to show that $L$ vanishes at these points. 
In the first case, $x_i=x_j$,  $j\neq i$, this is clear since
 the $i$th and $j$th rows in the determinant are equal.
It then follows from Lemma \ref{fl} that 
$L$ vanishes also at $x_i=1/tx_1\dotsm 
\hat x_i\dotsm x_n$. 

In the other cases, the same proof works with  obvious
modifications. It is actually enough to go through this 
for  $R=D_n$, since the remaining five cases 
can then be deduced using Lemma \ref{fl}.
\end{proof}
In the case $R=D_n$, one may well attribute Proposition \ref{wp} to
Warnaar.  Although he only states it in a special case, see Proposition
\ref{bcdet} below,  his proof 
extends \emph{verbatim} to  the general case.

\begin{Remark}
Replacing $x_i$ by $x_i/{\sqrt[n]{t}}$ one sees that \eqref{awpi} is
  equivalent to its special case $t=1$. Thus, if we would redefine
    $W_{A_{n-1}}$  as $\theta(x_1\dotsm
  x_n)W_{A_{n-1}}(x)$, we could give a unified statement of
  Proposition~\ref{wp} for all root systems. We have chosen to
  formulate the result using the superfluous parameter $t$ since this
  seems convenient for applications, in particular to multidimensional
  hypergeometric series. 
\end{Remark}

\section{Elliptic determinant evaluations}\label{secell}

We do not consider Proposition \ref{wp}  a determinant
\emph{evaluation}, since we do not have a simple formula for the
constant $C$. From our perspective, the main use of  Proposition
\ref{wp} is to systematize our knowledge of elliptic determinant
evaluations, as corresponding to various special cases when this
constant can be computed.

\subsection{Warnaar's type $D$ determinant}

For comparison and completeness, we first review the following
determinant evaluation due to Warnaar~\cite[Lemma~5.3]{W}. 
Warnaar used it to obtain a summation formula for a 
 multidimensional elliptic hyper\-geometric
series; further related applications may be found in \cite{Ro1},
\cite{Ro}, \cite{RS1}, \cite{Sp}. In the limit $p\rightarrow 0$
 it reduces to Krattenthaler's
determinant  \cite[Lemma~34]{K1}, which has been a powerful
tool in the enumeration of, and computation of generating functions
for, restricted families of plane partitions and tableaux, see the
discussion of Lemmas 3--5 and Theorems 26--31 in \cite{K2}.  

Warnaar's determinant corresponds to the case of Proposition
  \ref{wp} when $R=D_n$ and
$$f_j(x)=P_{j}(x)\prod_{k=j+1}^n\theta(a_kx^\pm),$$
 with $P_{j}$  a $D_j$ theta function. Then, for
$x_i=a_i$, the matrix in \eqref{wpi} is triangular, so that its
 determinant, and thus the
constant $C$, can be computed. This leads to the following result.

\begin{Proposition}[(Warnaar) A $D$ type determinant evaluation]
\label{bcdet}
Let $x_1,\dots,x_n$ and $a_1,\dots,a_n$ be indeterminates.
For each $j=1,\dots,n$, let $P_j$ be a $D_j$ theta function. 
Then there holds
$$
\det_{1\leq i,j\leq n}\left(P_{j}(x_i)\prod_{k=j+1}^n\theta(a_kx_i^\pm)
\right)=\prod_{i=1}^n P_{i}(a_i)
\,\prod_{1\leq i<j\leq
  n}a_jx_j^{-1}\theta(x_jx_i^\pm).
$$
\end{Proposition}

The parameter $a_1$  is  introduced for 
convenience, its value
being immaterial since $P_1$ is constant. Similar remarks can be made
about many of our results below.

\begin{Corollary}[A $D$ type Cauchy determinant]
\label{frobc}
Let $x_1,\dots,x_n$ and $a_1,\dots,a_n$ be indeterminates.
Then there holds
$$
\det_{1\leq i,j\leq n}\left(\frac 1{\theta(a_jx_i^\pm)}\right)
=\frac{\prod_{1\leq i<j\leq n}
a_jx_j^{-1}\,\theta(x_jx_i^\pm,a_ia_j^\pm)}
{\prod_{i,j=1}^n\theta(a_jx_i^\pm)}.
$$
\end{Corollary}

\begin{proof}
Let $P_{j}(x)=\prod_{k=1}^{j-1}\ta(a_kx^\pm)$ in Proposition \ref{bcdet},  pull
$\prod_{k=1}^n\ta(a_kx_i^\pm)$ out of the $i$th row of the determinant
($i=1,\dots,n$) and divide both sides by $\prod_{i,j=1}^n\theta(a_jx_i^\pm)$.
\end{proof}

Corollary~\ref{frobc} was used by Rains~\cite{R}, \cite{R2} to obtain
transformations and recurrences for multiple elliptic hypergeometric
integrals. Perhaps surprisingly, it is equivalent to the classical
Cauchy determinant
$$\det_{1\leq i,j\leq n}\left(\frac 1{u_i+v_j}\right)
=\frac{\prod_{1\leq i<j\leq n}
(u_j-u_i)(v_j-v_i)}
{\prod_{i,j=1}^n(u_i+v_j)},$$
 see
 \cite{R2}.

Another simple consequence of Proposition~\ref{bcdet} is the 
 following  determinant evaluation, which is included
 here for possible future reference. Two related determinant evaluations,
  corresponding to the type $A$ root system
 and restricted to the polynomial case, were applied in \cite{S1} and \cite{S3}
to obtain multidimensional matrix inversions that played a major
role in the derivation of new summation formulae for multidimensional
basic hypergeometric series, see
 Remark~\ref{remadetcorr}. Eventually, 
Corollary~\ref{bcdetcor}  may have similar applications in the elliptic setting.

\begin{Corollary}[A $D$ type determinant evaluation]\label{bcdetcor}
Let $x_1,\dots,x_n$, $a_1,\dots,a_{n+1}$ and $b$ be indeterminates.
For each $j=1,\dots,n+1$, let $P_j$ be a $D_j$ theta function. Then there holds
\begin{multline*}
P_{n+1}(b)\det_{1\le i,j\le n}\!\left(P_{j}(x_i)
\prod_{k=j+1}^{n+1}\ta(a_kx_i^\pm)
-\frac{P_{n+1}(x_i)}{P_{n+1}(b)}P_{j}(b)
\prod_{k=j+1}^{n+1}\ta(a_kb^\pm)\right)\\
=\prod_{i=1}^{n+1}P_{i}(a_i)
\prod_{1\le i<j\le n+1}a_jx_j^{-1}\,\ta(x_jx_i^\pm),
\end{multline*}
where $x_{n+1}=b$.
\end{Corollary}
\begin{proof}
We proceed similarly as in the proof of Lemma~A.1 of \cite{S1}.
In particular, we utilize
$\det\begin{pmatrix}M&\eta\\\xi&\gamma\end{pmatrix}=\gamma
\det\big(M-\gamma^{-1}\eta\xi\big)$ (which is a special case
of a formula due to Sylvester~\cite{S}) applied to
$M=\big(P_{j}(x_i)\prod_{k=j+1}^{n+1}\ta(a_kx_i^\pm)\big)$,
$\xi=\big(P_{j}(b)\prod_{k=j+1}^{n+1}\ta(a_kb^\pm)\big)$,
$\eta=(P_{n+1}(x_i))$, $\gamma=P_{n+1}(b)$, and then apply Proposition~\ref{bcdet}.
\end{proof}

\subsection{An $A$ type determinant}

If one tries to imitate  the proof of Proposition~\ref{bcdet},  using
Proposition \ref{wp} for $B_n$, $B_n^\vee$, $C_n$, $C_n^\vee$ or $BC_n$, rather
than $D_n$, one will find results that are equivalent to 
 Proposition \ref{bcdet} in view of Lemma~\ref{fl}. 
However, for the  root system $A_{n-1}$ one 
obtains the following  new elliptic extension of the
Vandermonde determinant \eqref{awd}, see Remark \ref{wdrem}.

\begin{Theorem}[An $A$ type determinant evaluation]\label{adet}
Let $x_1,\dots,x_n$, $a_1,\dots,a_n$, and $t$ be indeterminates.
For each $j=1,\dots,n$, let $P_j$ be an $A_{j-1}$ theta function of
 norm $ta_1\dotsm a_j$. Then there holds
\begin{equation}\label{adetid}
\det_{1\le i,j\le n}\left(P_j(x_i)
\prod_{k=j+1}^n\ta(a_kx_i)\right)
=\frac{\ta(ta_1\dotsm a_nx_1\dotsm x_n)}{\ta(t)}
\prod_{i=1}^nP_i(1/a_i)
\prod_{1\le i<j\le n}a_jx_j\,\ta(x_i/x_j).
\end{equation}
\end{Theorem}

\begin{proof}
By the $A_{n-1}$ case of Proposition \ref{wp}, with $t$ replaced by $ta_1\dotsm a_n$,
\eqref{adetid} holds up to a factor independent of
$x_i$. To compute this constant one may let $x_i=1/a_i$, in
which case the matrix on the left-hand side
is  triangular.
\end{proof}

By Lemma \ref{fl}, we may without loss of generality 
 assume that
\begin{equation}\label{pjf}P_j(x)=\theta(b_{1j}x)\dotsm\theta(b_{jj}x),
\end{equation}
where 
$b_{1j}\dotsm b_{jj}=ta_1\dotsm a_j$. On the right-hand side of
\eqref{adetid}, we then have
$P_1(1/a_1)/\theta(t)=1$. After 
replacing $t$ by $t/a_1\dotsm a_n$, this gives
 the following equi\-valent form of 
 Theorem \ref{adet}:
$$
\det_{1\le i,j\le n}\left(\prod_{k=1}^j\theta(b_{kj}x_i)
\prod_{k=j+1}^n\ta(a_kx_i)\right)
=\ta(tx_1\dots x_n)
\prod_{i=2}^n\prod_{k=1}^i\theta(b_{ki}/a_i)
\prod_{1\le i<j\le n}a_jx_j\,\ta(x_i/x_j),
$$
where
$$b_{1j}\dotsm b_{jj}a_{j+1}\dotsm a_n=t,\qquad j=1,\dots,n. $$

If we make the further specialization 
$$(b_{1j},\dots,b_{jj})=(c_1,\dots,c_{j-1},b_j) $$
and then interchange $a_j$ and $c_j$,
we recover the following determinant evaluation due to Tarasov and
Varchenko. In a special case, it was also obtained by
Hasegawa~\cite[Lemma~1]{H}, who used it to compute the
trace of elliptic $L$-operators, leading to the elliptic
Ruijsenaars(--Macdonald) commuting difference operators, see
 \cite{Ru}.

\begin{Corollary}[Tarasov and Varchenko]\label{tvcor}
Let $x_1,\dots, x_n$, $a_1,\dots, a_{n-1}$, $b_1,\dots, b_n$,
$c_2,\dots, c_{n}$ and $t$
 be indeterminates, such that
$$a_1\dotsm a_{j-1}b_j c_{j+1}\dotsm c_n=t,\qquad j=1,\dots,n. $$ 
Then there holds
$$
\det_{1\le i,j\le
  n}\left(\prod_{k=1}^{j-1}\theta(a_kx_i)\cdot\theta(b_jx_i)
\prod_{k=j+1}^n\ta(c_kx_i)\right)
=\ta(tx_1\dotsm x_n)
\prod_{i=2}^n\theta(b_i/c_i)
\prod_{1\le i<j\le n}c_jx_j\,\ta(x_i/x_j,a_i/c_j).
$$
\end{Corollary}

Note that $\prod_{i=2}^n\theta(b_i/c_i)=\prod_{i=1}^{n-1}\theta(b_{i}/a_{i})$.

\begin{Remark}\label{tvrem}
Corollary \ref{tvcor} appears rather implicitly in
 \cite[Appendix B]{TV},  as a special case of a much more general result. 
More precisely, it is the case
 $\ell=1$  of an infinite family of evaluations for the determinants
\begin{equation}\label{tvid}
\det_{\mathfrak l,\mathfrak m\in\mathcal Z_\ell^n}
(J_\mathfrak l(u\rhd \mathfrak m)),\end{equation}
where rows and columns are labelled by the compositions
$$\mathcal Z_\ell^n=\{\mathfrak l=(\mathfrak l_1,\dots,\mathfrak l_n);
\,\mathfrak l_i\geq 0,\
\sum \mathfrak l_i=\ell\}.$$
When $\ell=1$, $\mathcal Z_\ell^n$ can be identified  with $\{1,\dots,n\}$
and one gets a ``usual'' determinant.
For an explanation of the other symbols in \eqref{tvid}, the
reader is kindly referred to \cite{TV}.
\end{Remark}

If we let $a_j=c_j$ in Corollary \ref{tvcor} and replace $t$ by
$ta_1\dotsm a_n$, so that $b_j=ta_j$, we recover the following 
 determinant evaluation  due to
Frobenius~\cite{Fr}.
 This identity  has found applications to Ruijsenaars
operators  \cite{Ru}, to multidimensional elliptic hypergeometric
series and integrals 
\cite{KN}, \cite{R} and to number theory \cite{Ro2}. 
It is closely related to the  denominator formula for certain affine
superalgebras, see \cite{Ro2}.
For a generalization
to higher genus Riemann surfaces, see  \cite[Corollary~2.19]{F}.

\begin{Corollary}[(Frobenius) An $A$ type Cauchy determinant evaluation]
\label{froa}
Let $x_1,\dots,x_n$, $a_1,\dots,a_n$ and $t$ be indeterminates.
Then there holds
$$
\det_{1\leq i,j\leq n}\left(\frac{\ta(ta_jx_i)}{\theta(t,a_jx_i)}\right)
=\frac{\ta(ta_1\dotsm a_nx_1\dotsm x_n)}{\ta(t)}
\frac{\prod_{1\leq i<j\leq n}
a_jx_j\,\theta(a_i/a_j,x_i/x_j)}
{\prod_{i,j=1}^n\theta(a_jx_i)}.
$$
\end{Corollary}

Finally, the following result is included
here for similar reasons as Corollary
\ref{bcdetcor}.

\begin{Corollary}[An $A$ type determinant evaluation]\label{adetcor}
Let $x_1,\dots,x_n$, $a_1,\dots,a_{n+1}$ and $b$ be indeterminates.
For each $j=1,\dots,n+1$, let $P_j$ be an $A_{j-1}$ theta function  of 
 norm $ta_1\dotsm a_j$. Then there holds
\begin{multline}\label{adetcorid}
P_{n+1}(b)\;\det_{1\le i,j\le n}\left(P_j(x_i)
\prod_{k=j+1}^{n+1}\ta(a_kx_i)-\frac{P_{n+1}(x_i)}{P_{n+1}(b)}
P_j(b)\prod_{k=j+1}^{n+1}\ta(a_kb)\right)\\
=\frac{\ta(tba_1\dotsm a_{n+1}x_1\dotsm x_n)}{\ta(t)}
\prod_{i=1}^{n+1}P_i(1/a_i)
\prod_{1\le i<j\le n+1}a_jx_j\,\ta(x_i/x_j),
\end{multline}
where $x_{n+1}=b$.
\end{Corollary}
\begin{proof}
Proceed as in the proof of Corollary~\ref{bcdetcor}
but apply Theorem~\ref{adet} instead of Proposition~\ref{bcdet}.
\end{proof}

\subsection{A  $C$ type determinant}

The following identity, associated to the affine root system of type
$C$, provides a new elliptic extension of the Weyl denominator formulas
\eqref{bwd}, \eqref{cwd} and \eqref{dwd}, see Remark \ref{wdrem}.

\begin{Theorem}[A $C$ type determinant evaluation]\label{cdet}
Let $x_1,\dots,x_n$, $a_1,\dots,a_n$, and $c_1,\dots, c_{n+2}$ be
indeterminates. For each $j=1,\dots,n$, let $P_j$ be an $A_{j-1}$
theta function  of norm
$$(c_1\dotsm c_{n+2}a_{j+1}\dotsm a_n)^{-1}.$$ Then there holds
\begin{multline}\label{dv}
\det_{1\leq i,j\leq n}\left(x_i^{-n-1}
\prod_{k=1}^{n+2}\ta(c_kx_i)\,
P_j(x_i)\prod_{k=j+1}^n\ta(a_kx_i)
-x_i^{n+1}\prod_{k=1}^{n+2}\ta(c_kx_i^{-1})\,
P_j(x_i^{-1})\prod_{k=j+1}^n\ta(a_kx_i^{-1})\right)\\
=\frac{a_1\dotsm a_n}
{\ta(c_1\dotsm c_{n+2}a_1\dotsm a_n)}
\prod_{i=1}^nP_i(1/a_i)
\prod_{1\leq i<j\leq n+2}\ta(c_ic_j)\prod_{i=1}^nx_i^{-1}\ta(x_i^2)
\prod_{1\le i<j\le n}a_jx_i^{-1}\,\ta(x_ix_j^{\pm}).
\end{multline}
\end{Theorem}

Equivalently, factoring $P_j$ as in \eqref{pjf}, we have
\begin{multline*}
\det_{1\leq i,j\leq n}\left(x_i^{-n-1}
\prod_{k=1}^{n+2}\ta(c_kx_i)\,
\prod_{k=1}^j\theta(b_{kj}x_i)\prod_{k=j+1}^n\ta(a_kx_i)\right.\\
\left.-x_i^{n+1}\prod_{k=1}^{n+2}\ta(c_kx_i^{-1})\,
\prod_{k=1}^j\theta(b_{kj}x_i^{-1})\prod_{k=j+1}^n\ta(a_kx_i^{-1})\right)\\
=-\frac{1}
{ c_1\dotsm c_{n+2}}
\prod_{i=2}^n\prod_{k=1}^i\theta(b_{ki}/a_i)
\prod_{1\leq i<j\leq n+2}\ta(c_ic_j)\prod_{i=1}^n x_i^{-1}\ta(x_i^2)
\prod_{1\le i<j\le n}a_jx_i^{-1}\,\ta(x_ix_j^{\pm}),
\end{multline*}
where 
$$b_{1j}\dotsm b_{jj}a_{j+1}\dotsm a_nc_1\dotsm c_{n+2}=1, 
\qquad j=1,\dots,n.$$

We will give two proofs of Theorem \ref{cdet}.

\begin{proof}[First proof of Theorem \ref{cdet}]
Using Lemma \ref{dl}, one  checks that the determinant  is of the form
\eqref{wpi}, with $R=C_n$.  Proposition \ref{wp} then guarantees that
the quotient of the two sides of \eqref{dv} is a constant, so it is
enough to verify the equality for some fixed values of $x_i$. We
choose $x_i=c_i$, so that  the second term in each matrix element
vanishes. The factor
$\prod_{k=1}^{n+2}\theta(c_kx_i)$ may then be pulled out from the $i$th row of the
determinant and cancelled, using
 $$\prod_{i=1}^n\prod_{k=1}^{n+2}\theta(c_kx_i)=
\frac 1{\theta(c_{n+1}c_{n+2})}
\prod_{1\leq i<j\leq n+2}\ta(c_ic_j)\prod_{i=1}^n\ta(x_i^2)\prod_{1\le
  i<j\le n}\ta(x_ix_j).$$
Introducing the parameter  $t=1/c_1\dotsm c_{n+2}a_1\dotsm
a_n$, we note that
$$\frac{\theta(c_{n+1}c_{n+2})}{\theta(c_1\dotsm c_{n+2}a_1\dotsm
a_n)}=\frac{\ta(ta_1\dotsm a_nx_1\dotsm x_n)}{\ta(t)}\prod_{i=1}^n\frac{1}{a_ix_i}. $$
Thus, we are reduced to proving
$$
\det_{1\leq i,j\leq n}\left(P_j(x_i)\prod_{k=j+1}^n\ta(a_kx_i)\right)
=\frac{\ta(ta_1\dotsm a_nx_1\dotsm x_n)}
{\ta(t)}
\prod_{i=1}^nP_i(1/a_i)
\prod_{1\le i<j\le n}a_jx_j\,\ta(x_i/x_j),
$$
where  $P_j$ is an $A_{j-1}$ theta function of norm $ta_1\dotsm
a_j$, and where $x_j$ may again be viewed as free variables.  This
 is exactly Theorem \ref{adet}.
\end{proof}

Let $\mathcal R_i$  denote the reflection operator $\mathcal R_if(x_i)=f(x_i^{-1})$. 
Then, due to linearity of the determinant, the  left-hand side of \eqref{dv} may be written
\begin{multline}\label{hcd}\prod_{i=1}^n(1-\mathcal R_i)\prod_{i=1}^n\left(x_i^{-n-1}\prod_{k=1}^{n+2}\ta(c_kx_i)
\right)\det_{1\leq i,j\leq n}\left(P_j(x_i)\prod_{k=j+1}^n\ta(a_kx_i)
\right)\\
=\frac{1}{\theta(1/c_1\dotsm c_{n+2}a_1\dotsm
a_n)}\prod_{i=1}^nP_i(1/a_i)\\
\times
\prod_{i=1}^n(1-\mathcal R_i)\,\theta\left(\frac{x_1\dotsm x_n}{c_1\dotsm c_{n+2}}\right)
\prod_{i=1}^n\left(x_i^{-n-1}\prod_{k=1}^{n+2}\ta(c_kx_i)
\right)\prod_{1\le i<j\le n}a_jx_j\,\ta(x_i/x_j),
 \end{multline} 
where we used Theorem \ref{adet} to compute the determinant.
Comparing this with the right-hand side of \eqref{dv} gives
the following equivalent form of Theorem \ref{cdet}.

\begin{Corollary}\label{lem}
In the notation above,
\begin{multline*}
\prod_{i=1}^n(1-\mathcal R_i)\,
\ta\!\left(\frac{x_1\dotsm x_n}{c_1\dotsm c_{n+2}}\right)
\prod_{i=1}^n\left( x_i^{-n-1}\prod_{j=1}^{n+2}\ta(c_jx_i)\right)
\prod_{1\le i<j\le n}x_j\,\ta(x_i/x_j)\\
=-\frac{1}{c_1\dotsm c_{n+2}}
\prod_{1\leq i<j\leq n+2}\ta(c_ic_j)\prod_{i=1}^n x_i^{-1}\ta(x_i^2)
\prod_{1\leq i<j\leq n}
x_i^{-1}\,\ta(x_ix_j^\pm).
\end{multline*}
\end{Corollary}

Corollary \ref{lem} resembles some  identities in the work of
Rains~\cite{R}.
It can be used to give an alternative proof of his type I $BC_n$
integral, originally conjectured by van Diejen and Spiridonov \cite{DS} 
 (Rains, personal communication).
It would be interesting to know if Corollary \ref{lem} can  be obtained by
specializing  a multidimensional elliptic hypergeometric summation
theorem on $0\leq k_i\leq m_i$ ($i=1,\dots,n$) to the case
$m_i\equiv 1$.

One consequence of  \eqref{hcd} is that if we can compute the
 left-hand side  for  some special choice of  $a_j$
and $P_j$, we can  compute it in general, since $a_j$ and $P_j$ 
 appear trivially on the right-hand side.
This observation can be used to give 
an alternative proof of
Theorem \ref{cdet}, based on the type $D$ Cauchy determinant of
Corollary~\ref{frobc}.

\begin{proof}[Second proof of Theorem~\ref{cdet}]
 We  consider the special case when $a_j=c_j^{-1}$, $1\leq j\leq n$, and
$$P_j(x)=\ta(tc_j^{-1}x)\prod_{k=1}^{j-1}\ta(c_k^{-1}x),$$
where $tc_{n+1}c_{n+2}=1$.
Then, the left-hand side of \eqref{dv} can be written
$$\det_{1\leq i,j\leq n}
\left(\left(\prod_{k=1,\,k\neq j}^n c_k^{-1}\,\ta(c_kx_i^\pm)\right)
(1-\mathcal R_i)\,x_i^{-2}\,\ta(c_{n+1}x_i,c_{n+2}x_i,c_jx_i,tc_j^{-1}x_i)
\right).$$
By  \eqref{radd} and 
Corollary~\ref{frobc}, this equals
\begin{multline*}\det_{1\leq i,j\leq n}
\left(\left(\prod_{k=1,\,k\neq j}^n c_k^{-1}\,\ta(c_kx_i^\pm)\right)
x_i^{-1}c_j^{-1}\,\ta(x_i^2,t,c_jc_{n+1},c_jc_{n+2})
\right)\\
=\frac{\ta(t)^n}{c_1^n\dotsm c_n^n}\prod_{i=1}^n x_i^{-1}\,
\ta(x_i^2,c_ic_{n+1},c_ic_{n+2})\prod_{i,j=1}^n\ta(c_jx_i^\pm)
\det_{1\leq i,j\leq n}\left(\frac 1{\ta(c_jx_i^\pm)}
\right)\\
=\frac{\ta(t)^n}{c_1^n\dotsm c_n^n}\prod_{i=1}^n x_i^{-1}\,
\ta(x_i^2,c_ic_{n+1},c_ic_{n+2})
\prod_{1\leq i<j\leq n}
c_jx_j^{-1}\,\ta(x_jx_i^\pm,c_ic_j^\pm),
\end{multline*}
which agrees with the right-hand side of \eqref{dv}.
As was remarked above, the general case now follows using \eqref{hcd}.
\end{proof}

\subsection{Determinants of type
 $B$,  $B^\vee$, $C^\vee$, $BC$ and $D$} 

If $c^2\in p^{\mathbb Z}$, then 
$\theta(cx)$ and $\theta(c/x)$ are equal up to a trivial factor.
Thus, if  one of the parameters $c_j$ in Theorem \ref{cdet}
is of this form, then the factor $\prod_{i=1}^n\theta(c_jx_i)$ may be pulled
out from the determinant. Up to the trivial scaling $c_j\mapsto pc_j$,
there are four choices: $c_j\in\{1,-1, p^{\frac
  12},-p^{\frac 12}\}$. By \eqref{txsq}, $\theta(c_jx_i)$ then
cancels against a part of the factor $\theta(x_i^2)$ on the right-hand
side. Making various specializations of this sort, the 
  $C_n$ Macdonald denominator in \eqref{dv} can be reduced to the 
 Macdonald denominator for  $B_n$,  $B_n^\vee$, $C_n^\vee$, $BC_n$
 and $D_n$. 

As a first example, we let  $c_{n+2}=-1$ in Theorem \ref{cdet}. Then,
$$\frac{\theta(x_i^2)}{\theta(c_{n+2}x_i)}=\theta(x_i)\theta(px_i^2;p^2).$$
 This gives the following determinant of type $BC$.

\begin{Corollary}[A $BC$ type determinant evaluation]\label{bcdet2}
Let $x_1,\dots,x_n$, $a_1,\dots,a_n$, and $c_1,\dots, c_{n+1}$ be
indeterminates. For each $j=1,\dots,n$, let $P_j$ be an $A_{j-1}$
theta function  of  norm
$$-(c_1\dotsm c_{n+1}a_{j+1}\dotsm a_n)^{-1}.$$
 Then there holds
\begin{multline*}
\det_{1\leq i,j\leq n}\left(x_i^{-n}
\prod_{k=1}^{n+1}\ta(c_kx_i)\,
P_j(x_i)\prod_{k=j+1}^n\ta(a_kx_i)
-x_i^{n+1}\prod_{k=1}^{n+1}\ta(c_kx_i^{-1})\,
P_j(x_i^{-1})\prod_{k=j+1}^n\ta(a_kx_i^{-1})\right)\\
=\frac{a_1\dotsm a_n}
{\ta(-c_1\dotsm c_{n+1}a_1\dotsm a_n)}
\prod_{i=1}^nP_i(1/a_i)\\\times
\prod_{i=1}^{n+1}\theta(-c_i)
\prod_{1\leq i<j\leq n+1}\ta(c_ic_j)\prod_{i=1}^n\ta(x_i)\ta(px_i^2;p^2)
\prod_{1\le i<j\le n}a_jx_i^{-1}\,\ta(x_ix_j^{\pm}).
\end{multline*}
\end{Corollary}

If we let $c_{n+1}=-p^{\frac12}$ in Corollary \ref{bcdet2},  we obtain
the following determinant of type $C^\vee$.

\begin{Corollary}[A $C^\vee$ type determinant evaluation]\label{cvdet}
Let $x_1,\dots,x_n$, $a_1,\dots,a_n$, and $c_1,\dots, c_{n}$ be
indeterminates. For each $j=1,\dots,n$, let $P_j$ be an $A_{j-1}$
theta function of norm
$$(p^{\frac 12} c_1\dotsm c_{n}a_{j+1}\dotsm a_n)^{-1}.$$
 Then there holds
\begin{multline*}
\det_{1\leq i,j\leq n}\left(x_i^{-n}
\prod_{k=1}^{n}\ta(c_kx_i)\,
P_j(x_i)\prod_{k=j+1}^n\ta(a_kx_i)
-x_i^{n+1}\prod_{k=1}^{n}\ta(c_kx_i^{-1})\,
P_j(x_i^{-1})\prod_{k=j+1}^n\ta(a_kx_i^{-1})\right)\\
=\frac{a_1\dotsm a_n\theta(p^{\frac 12})}
{\ta(p^{\frac 12}c_1\dotsm c_{n}a_1\dotsm a_n)}
\prod_{i=1}^nP_i(1/a_i)\\
\times\prod_{i=1}^n\ta(-c_i,p^\frac12)
\prod_{1\leq i<j\leq n}\ta(c_ic_j)\prod_{i=1}^n\ta(x_i;p^{\frac 12})
\prod_{1\le i<j\le n}a_jx_i^{-1}\,\ta(x_ix_j^{\pm}).
\end{multline*}
\end{Corollary}

If we let  $c_{n+1}=- p^\frac12$ and
$c_{n+2}=p^\frac12$ in Theorem \ref{cdet}, and replace $c_1$ by
$c_1/p$ for convenience,  we obtain the following determinant of type
$B^\vee$.

\begin{Corollary}[A $B^\vee$ type determinant evaluation]\label{bvdet}
Let $x_1,\dots,x_n$, $a_1,\dots,a_n$, and $c_1,\dots, c_{n}$ be
indeterminates. For each $j=1,\dots,n$, let $P_j$ be an $A_{j-1}$
theta function of norm
$$-(c_1\dotsm c_{n}a_{j+1}\dotsm a_n)^{-1}.$$
 Then there holds
\begin{multline*}
\det_{1\leq i,j\leq n}\left(x_i^{-n}
\prod_{k=1}^{n}\ta(c_kx_i)\,
P_j(x_i)\prod_{k=j+1}^n\ta(a_kx_i)
-x_i^{n}\prod_{k=1}^{n}\ta(c_kx_i^{-1})\,
P_j(x_i^{-1})\prod_{k=j+1}^n\ta(a_kx_i^{-1})\right)\\
=\frac{a_1\dotsm a_nc_1\dotsm c_n\theta(-1)}
{\ta(-c_1\dotsm c_{n}a_1\dotsm a_n)}
\prod_{i=1}^nP_i(1/a_i)\\\times
\prod_{i=1}^n \ta(pc_i^2;p^2)\prod_{1\leq i<j\leq n}\ta(c_ic_j)\prod_{i=1}^n x_i^{-1}\ta(x_i^2;p^2)
\prod_{1\le i<j\le n}a_jx_i^{-1}\,\ta(x_ix_j^{\pm}).
\end{multline*}
\end{Corollary}

If we let 
$c_n=-1$ in Corollary \ref{bvdet} we obtain, using also \eqref{tev}, 
the following determinant
of type $B$.

\begin{Corollary}[A $B$ type determinant evaluation]\label{bdet}
Let $x_1,\dots,x_n$, $a_1,\dots,a_n$, and $c_1,\dots, c_{n-1}$ be
indeterminates. For each $j=1,\dots,n$, let $P_j$ be an $A_{j-1}$
theta function  of  norm
$$(c_1\dotsm c_{n-1}a_{j+1}\dotsm a_n)^{-1}.$$
 Then there holds
\begin{multline*}
\det_{1\leq i,j\leq n}\left(x_i^{1-n}
\prod_{k=1}^{n-1}\ta(c_kx_i)\,
P_j(x_i)\prod_{k=j+1}^n\ta(a_kx_i)
-x_i^{n}\prod_{k=1}^{n-1}\ta(c_kx_i^{-1})\,
P_j(x_i^{-1})\prod_{k=j+1}^n\ta(a_kx_i^{-1})\right)\\
=-\frac{2a_1\dotsm a_nc_1\dotsm c_{n-1}}
{\ta(c_1\dotsm c_{n-1}a_1\dotsm a_n)}
\prod_{i=1}^nP_i(1/a_i)\\\times
\prod_{i=1}^{n-1}\theta(-c_i)
\ta(pc_i^2;p^2)\prod_{1\leq i<j\leq n-1}\ta(c_ic_j)
\prod_{i=1}^n\theta(x_i)
\prod_{1\le i<j\le n}a_jx_i^{-1}\,\ta(x_ix_j^{\pm}).
\end{multline*}
\end{Corollary}

Finally,  assuming $n\geq 2$, we  let $c_{n-1}=1$ in Corollary
\ref{bdet}. 
Again using \eqref{tev}, we  obtain 
following type $D$ determinant.

\begin{Corollary}[A $D$ type determinant evaluation]\label{ddet}
Let $x_1,\dots,x_n$, $a_1,\dots,a_n$, and $c_1,\dots, c_{n-2}$ be
indeterminates. For each $j=1,\dots,n$, let $P_j$ be an $A_{j-1}$
theta function  of norm
$$(c_1\dotsm c_{n-2}a_{j+1}\dotsm a_n)^{-1}.$$
 Then, for $n\geq 2$,
there holds
\begin{multline*}
\det_{1\leq i,j\leq n}\left(x_i^{1-n}
\prod_{k=1}^{n-2}\ta(c_kx_i)\,
P_j(x_i)\prod_{k=j+1}^n\ta(a_kx_i)
+x_i^{n-1}\prod_{k=1}^{n-2}\ta(c_kx_i^{-1})\,
P_j(x_i^{-1})\prod_{k=j+1}^n\ta(a_kx_i^{-1})\right)\\
=-\frac{4a_1\dotsm a_nc_1\dotsm c_{n-2}}
{\ta(c_1\dotsm c_{n-2}a_1\dotsm a_n)}
\prod_{i=1}^nP_i(1/a_i)
\prod_{1\leq i\leq j\leq n-2}\ta(c_ic_j)
\prod_{1\le i<j\le n}a_jx_i^{-1}\,\ta(x_ix_j^{\pm}).
\end{multline*}
\end{Corollary}

\section{Some polynomial determinant evaluations}\label{secrat}

In this Section we consider the polynomial special case, $p=0$, of the elliptic
determinant evaluations in Section \ref{secell}. The resulting
identities involve
the Weyl denominator of classical (non-affine) root systems, cf.\ \eqref{wd}.

We must first interpret the term ``$A_{n-1}$  theta function''
 in the case $p=0$. 
One way is to rewrite Definition \ref{rtd} 
in terms of the Laurent coefficients of
 $f(x)=\sum_j a_jx^j$. Namely, $f$ is an  $A_{n-1}$  theta function of
 norm $t$ if and only if
$$a_{j+n}=(-1)^ntp^ja_j. $$
When $p=0$ this means  that $a_j=0$ unless $0\leq j\leq n$ and that
$a_n=(-1)^nta_0$. Thus, we obtain precisely the space of
polynomials  of degree $n$ and \emph{norm} $t$, where the norm of
$a_0+a_1x+\dots+a_nx^n$ is 
defined as $(-1)^na_n/a_0$. Equivalently, the polynomial
 $C(1-b_1x)\dotsm (1-b_nx)$ has norm $b_1\dotsm
 b_n$. Thus, we obtain the same result by formally letting $p=0$ in
 Lemma \ref{fl}. 
With this interpretation of the term
$A_{n-1}$  theta function, Theorems \ref{adet} and \ref{cdet}
 remain valid when $p=0$. 

\subsection{Determinants of type $A$}

We first give 
the case $p=0$ of Theorem~\ref{adet}.

\begin{Corollary}[An $A$ type determinant evaluation]\label{apoldet}
Let $x_1,\dots,x_n$, $a_1,\dots,a_n$, and $t$ be indeterminates.
For each $j=1,\dots,n$, let $P_j$ be a polynomial
 of degree $j$ and norm $ta_1\dotsm a_j$. Then there holds
$$
\det_{1\le i,j\le n}\left(P_j(x_i)
\prod_{k=j+1}^n(1-a_kx_i)\right)
=\frac{1-ta_1\dotsm a_nx_1\dotsm x_n}{1-t}
\prod_{i=1}^nP_i(1/a_i)
\prod_{1\le i<j\le n}a_j(x_j-x_i).
$$
\end{Corollary}

 It is easy to prove Corollary \ref{apoldet}
directly by a standard ``identification of
factors'' argument.

It is possible to remove the restriction on the norm of the
polynomials $P_j$ through
a limit transition, decreasing their degree  by one. 
Such limits do not make sense in the elliptic case ($p\neq 0$).
This leads to the following
determinant evaluation due to Krattenthaler \cite[Lemma~35]{K1}, 
who obtained it as a limit case of  \cite[Lemma~34]{K1},
see the discussion of 
Proposition \ref{bcdet} above.

\begin{Corollary}[(Krattenthaler)\label{apoldet2}
An $A$ type determinant evaluation]
Let $x_1,\dots,x_n$ and $a_1,\dots,a_n$ be indeterminates.
For each $j=1,\dots,n$, let $P_{j-1}$ be a polynomial
 of degree at most $j-1$. Then there holds
$$
\det_{1\le i,j\le n}\left(P_{j-1}(x_i)
\prod_{k=j+1}^n(1-a_kx_i)\right)
=\prod_{i=1}^nP_{i-1}(1/a_i)
\prod_{1\le i<j\le n}a_j(x_j-x_i).
$$
\end{Corollary}

\begin{proof}
In Corollary \ref{apoldet}, write $P_j(x)=(1-tb_jx)\tilde P_{j-1}(x)$,
 let $t\rightarrow 0$ and then relabel  $\tilde P_{j-1}\mapsto P_{j-1}$.
\end{proof}

We also note the following consequence of Corollary \ref{adetcor}.

\begin{Corollary}[An $A$ type determinant evaluation]\label{adetcorr}
Let $x_1,\dots,x_n$ and $b$ be indeterminates.
For each $j=1,\dots,n$, let $P_{j-1}(x)$ be a polynomial  in
$x$ of degree at most $j-1$ with constant term $1$,
and let $Q(x)=(1-y_1x)\dotsm (1-y_{n+1}x)$.
Then there holds
\begin{multline}\label{adetcorrid}
Q(b)\;\det_{1\le i,j\le n}\left(x_i^{n+1-j}P_{j-1}(x_i)
-b^{n+1-j}P_{j-1}(b)\frac{Q(x_i)}{Q(b)}\right)\\
=(1-bx_1\dots x_ny_1\dotsm y_{n+1})
\prod_{i=1}^n(x_i-b)\prod_{1\le i<j\le n}(x_i-x_j).
\end{multline}
\end{Corollary}

\begin{proof}
In Corollary~\ref{adetcor},  let $p=0$ and assume, as a matter of
normalization, that the polynomials $P_j$ have constant term $1$.
Write $t=s^{n+1}$, $a_i=c_i/s$,
$$P_j(x)=(1-s^{n+1-j}d_jx)\tilde P_{j-1}(x),\qquad j=1,\dots, n,$$
$$P_{n+1}(x)=(1-y_1x)\dotsm (1-y_{n+1}x).$$
Then,   
 $\tilde P_{j-1}$ has norm $c_1\dotsm c_j/d_j$
and $P_{n+1}$ norm $y_1\dotsm y_{n+1}=c_1\dotsm c_{n+1}$,
 which are in particular
 independent of $s$. 
Dividing  both sides of \eqref{adetcorid} by 
$\prod_{1\le i<j\le n+1}(-a_j),$
 letting  $s\rightarrow 0$ and finally
relabelling $\tilde P_{j-1}\mapsto P_{j-1}$, $P_{n+1}\mapsto Q$,
 we obtain the desired result.
\end{proof}

\begin{Remark}\label{remadetcorr}
Note that the right-hand side of \eqref{adetcorrid} is independent of
$P_{j-1}$. The special case $P_{j-1}(x)=1$, for $j=1,\dots,n$,
is Lemma~A.1 of \cite{S1}, which was needed in order to obtain an
$A_n$ matrix inversion that played a crucial role in the derivation
of multiple basic hypergeometric series identities.
A slight generalization of \cite[Lemma~A.1]{S1} was given
in \cite[Lemma~A.1]{S3}.
\end{Remark}

\subsection{Determinants of type $B$, $C$, and $D$}

Next, we turn to the
 $p=0$ case of Theorem~\ref{cdet}.

\begin{Corollary}[A $C$ type determinant evaluation]\label{cdetr}
Let $x_1,\dots,x_n$, $a_1,\dots,a_n$, and $c_1,\dots,c_{n+2}$ be
indeterminates. For each $j=1,\dots,n$, let $P_j$ be a polynomial
  of degree $j$ with norm
$$(c_1\dotsm c_{n+2}a_{j+1}\dotsm a_n)^{-1}.$$
Then there holds
\begin{multline*}
\det_{1\leq i,j\leq n}\left(x_i^{-n-1}
\prod_{k=1}^{n+2}(1-c_kx_i)\,
P_j(x_i)\prod_{k=j+1}^n(1-a_kx_i)\right.\\
\left.-x_i^{n+1}\prod_{k=1}^{n+2}(1-c_kx_i^{-1})\,
P_j(x_i^{-1})\prod_{k=j+1}^n(1-a_kx_i^{-1})\right)\\
=\frac{a_1\dotsm a_n}
{1-c_1\dotsm c_{n+2}a_1\dotsm a_n}
\prod_{i=1}^nP_i(1/a_i)
\prod_{1\leq i<j\leq n+2}(1-c_ic_j)\prod_{i=1}^n x_i^{-n}(1-x_i^2)
\prod_{1\le i<j\le n}a_j(x_j-x_i)(1-x_ix_j).
\end{multline*}
\end{Corollary}

If we let $c_{n+2}=-1$ in Corollary \ref{cdetr} or, equivalently,
$p=0$ in Corollary \ref{bcdet2}, we obtain the following determinant
of type $B$. 

\begin{Corollary}[A $B$ type determinant evaluation]\label{bdetr}
Let $x_1,\dots,x_n$, $a_1,\dots,a_n$, and $c_1,\dots,c_{n+1}$ be
indeterminates. For each $j=1,\dots,n$, let $P_j$ be a polynomial
  of degree $j$ with norm
$$-(c_1\dotsm c_{n+1}a_{j+1}\dotsm a_n)^{-1}.$$
Then there holds
\begin{multline*}
\det_{1\leq i,j\leq n}\left(x_i^{-n}
\prod_{k=1}^{n+1}(1-c_kx_i)\,
P_j(x_i)\prod_{k=j+1}^n(1-a_kx_i)\right.\\
\left.-x_i^{n+1}\prod_{k=1}^{n+1}(1-c_kx_i^{-1})\,
P_j(x_i^{-1})\prod_{k=j+1}^n(1-a_kx_i^{-1})\right)\\
=\frac{a_1\dotsm a_n}
{1+c_1\dotsm c_{n+1}a_1\dotsm a_n}
\prod_{i=1}^nP_i(1/a_i)\\
\times\prod_{1\leq i<j\leq n+1}(1-c_ic_j)
\prod_{i=1}^{n+1}(1+c_i)
\prod_{i=1}^n x_i^{1-n}(1-x_i)
\prod_{1\le i<j\le n}a_j(x_j-x_i)(1-x_ix_j).
\end{multline*}
\end{Corollary}

If we let $c_{n+1}=1$ in Corollary \ref{bdetr}, the factor
$\prod_{i=1}^n(1-x_i)$ may be cancelled. This gives the following
determinant of type $D$.

\begin{Corollary}[A $D$ type determinant evaluation]\label{ddetr}
Let $x_1,\dots,x_n$, $a_1,\dots,a_n$, and $c_1,\dots,c_n$ be
indeterminates. For each $j=1,\dots,n$, let $P_j$ be a polynomial
 of degree $j$ with norm
$$-(c_1\dotsm c_na_{j+1}\dotsm a_n)^{-1}.$$
Then there holds
\begin{multline*}
\det_{1\leq i,j\leq n}\left(x_i^{-n}
\prod_{k=1}^n(1-c_kx_i)\,
P_j(x_i)\prod_{k=j+1}^n(1-a_kx_i)\right.\\
\left.+x_i^n\prod_{k=1}^n(1-c_kx_i^{-1})\,
P_j(x_i^{-1})\prod_{k=j+1}^n(1-a_kx_i^{-1})\right)\\
=\frac{2\,a_1\dotsm a_n}
{1+c_1\dotsm c_na_1\dotsm a_n}
\prod_{i=1}^n  P_i(1/a_i)
\prod_{1\leq i\le j\leq n}(1-c_ic_j)\prod_{i=1}^n x_i^{1-n}
\prod_{1\le i<j\le n}a_j(x_j-x_i)(1-x_ix_j).
\end{multline*}
\end{Corollary}

Similarly as when deriving Corollary \ref{apoldet2} from Corollary
\ref{apoldet}, we may  remove the  restriction on the norm of  $P_j$
 in Corollaries~\ref{cdetr}, \ref{bdetr} and \ref{ddetr} by a limit
 transition, through which their degree is lowered by one.

\begin{Corollary}[A $C$ type determinant evaluation]\label{cdetr1}
Let $x_1,\dots,x_n$, $a_1,\dots,a_n$, and $c_1,\dots,c_{n+1}$ be
indeterminates. For each $j=1,\dots,n$, let $P_{j-1}$ be a polynomial
  of degree at most $j-1$. Then there holds
\begin{multline*}
\det_{1\leq i,j\leq n}\left(x_i^{-n}
\prod_{k=1}^{n+1}(1-c_kx_i)\,
P_{j-1}(x_i)\prod_{k=j+1}^n(1-a_kx_i)\right.\\
\left.-x_i^n\prod_{k=1}^{n+1}(1-c_kx_i^{-1})\,
P_{j-1}(x_i^{-1})\prod_{k=j+1}^n(1-a_kx_i^{-1})\right)\\
=\prod_{i=1}^nP_{i-1}(1/a_i)\prod_{1\leq i<j\leq n+1}(1-c_ic_j)
\prod_{i=1}^nx_i^{-n}(1-x_i^2)
\prod_{1\le i<j\le n}a_j(x_j-x_i)(1-x_ix_j).
\end{multline*}
\end{Corollary}

\begin{proof}
In Corollary \ref{cdetr}, write $P_j(x)=(x+b_jc_{n+2})\tilde
P_{j-1}(x)$, let $c_{n+2}\rightarrow 0$ and relabel $\tilde
P_{j-1}\mapsto P_{j-1}$.
\end{proof}

\begin{Corollary}[A $B$ type determinant evaluation]\label{bdetr1}
Let $x_1,\dots,x_n$, $a_1,\dots,a_n$, and $c_1,\dots,c_n$ be
indeterminates. For each $j=1,\dots,n$, let $P_{j-1}$ be a polynomial
  of degree at most $j-1$. Then there holds
\begin{multline*}
\det_{1\leq i,j\leq n}\left(x_i^{1-n}
\prod_{k=1}^n(1-c_kx_i)\,
P_{j-1}(x_i)\prod_{k=j+1}^n(1-a_kx_i)\right.\\
\left.-x_i^n\prod_{k=1}^n(1-c_kx_i^{-1})\,
P_{j-1}(x_i^{-1})\prod_{k=j+1}^n(1-a_kx_i^{-1})\right)\\
=\prod_{i=1}^nP_{i-1}(1/a_i)
\prod_{1\leq i<j\leq n}(1-c_ic_j)\prod_{i=1}^n(1+c_i)
\prod_{i=1}^nx_i^{1-n}(1-x_i)
\prod_{1\le i<j\le n}a_j(x_j-x_i)(1-x_ix_j).
\end{multline*}
\end{Corollary}

\begin{proof}
Let $c_{n+1}=-1$ in Corollary \ref{cdetr1} and divide by
$\prod_{i=1}^n(1+x_i^{-1})$. 
\end{proof}

\begin{Corollary}[A $D$ type determinant evaluation]\label{ddetr1}
Let $x_1,\dots,x_n$, $a_1,\dots,a_n$, and $c_1,\dots,c_{n-1}$ be
indeterminates. For each $j=1,\dots,n$, let $P_{j-1}$ be a polynomial
  of degree at most $j-1$. Then there holds
\begin{multline*}
\det_{1\leq i,j\leq n}\left(x_i^{1-n}
\prod_{k=1}^{n-1}(1-c_kx_i)\,
P_{j-1}(x_i)\prod_{k=j+1}^n(1-a_kx_i)\right.\\
\left.+x_i^{n-1}\prod_{k=1}^{n-1}(1-c_kx_i^{-1})\,
P_{j-1}(x_i^{-1})\prod_{k=j+1}^n(1-a_kx_i^{-1})\right)\\
=2\prod_{i=1}^n P_{i-1}(1/a_i)\prod_{1\leq i\le j\leq n-1}(1-c_ic_j)
\prod_{i=1}^n x_i^{1-n}\prod_{1\le i<j\le n}a_j(x_j-x_i)(1-x_ix_j).
\end{multline*}
\end{Corollary}

\begin{proof}
Let $c_n=1$ in Corollary \ref{bdetr1} and divide by $\prod_{i=1}^n(1-x_i)$.
\end{proof}

Next, we give  some further specializations of our determinant
evaluations, which
are closer to the classical Weyl denominator formulas.

\begin{Corollary}[A $C$ type determinant evaluation]\label{cdetr1cor}
Let $x_1,\dots,x_n$, and $c_1,\dots,c_{n+1}$ be
indeterminates. For each $j=1,\dots,n$, let $P_{j-1}$ be a polynomial
 of degree at most $j-1$. Then there holds
\begin{multline*}
\det_{1\leq i,j\leq n}\left(x_i^{-j}
\prod_{k=1}^{n+1}(1-c_kx_i)\,P_{j-1}(x_i)
-x_i^j\prod_{k=1}^{n+1}(1-c_kx_i^{-1})\,P_{j-1}(x_i^{-1})\right)\\
=\prod_{i=1}^n P_{i-1}(0)
\prod_{1\leq i<j\leq n+1}(1-c_ic_j)
\prod_{i=1}^nx_i^{-n}(1-x_i^2)
\prod_{1\le i<j\le n}(x_i-x_j)(1-x_ix_j).
\end{multline*}
\end{Corollary}
\begin{proof}
In Corollary~\ref{cdetr1},
divide both sides of the identity by $\prod_{1\le i<j\le n}(-a_j)$,
and then let $a_j\to\infty$, successively for $j=2,\dots,n$.
\end{proof}

\begin{Remark} 
The special case $P_{j-1}(x)=1$, for $j=1,\dots,n$,
is Lemma~A.11 of \cite{S1}, needed in order to obtain a $C_n$ matrix
inversion (which was later applied in \cite{S2}).
\end{Remark}

\begin{Corollary}[A $B$ type determinant evaluation]\label{bdetr1cor}
Let $x_1,\dots,x_n$ and $c_1,\dots,c_n$ be
indeterminates. For each $j=1,\dots,n$, let $P_{j-1}$ be a polynomial
  of degree at most $j-1$. Then there holds
\begin{multline*}
\det_{1\leq i,j\leq n}\left(x_i^{1-j}
\prod_{k=1}^n(1-c_kx_i)\,P_{j-1}(x_i)
-x_i^j\prod_{k=1}^n(1-c_kx_i^{-1})\,P_{j-1}(x_i^{-1})\right)\\
=\prod_{i=1}^n P_{i-1}(0)\prod_{1\leq i<j\leq n}(1-c_ic_j)\prod_{i=1}^n(1+c_i)
\prod_{i=1}^n x_i^{1-n}(1-x_i)
\prod_{1\le i<j\le n}(x_i-x_j)(1-x_ix_j).
\end{multline*}
\end{Corollary}

\begin{proof}
Let $c_{n+1}=-1$ in Corollary \ref{cdetr1cor} and divide by $\prod_{i=1}^n(1+x_i^{-1})$.
\end{proof}

\begin{Corollary}[A $D$ type determinant evaluation]\label{ddetr1cor}
Let $x_1,\dots,x_n$ and $c_1,\dots,c_{n-1}$ be indeterminates.
For each $j=1,\dots,n$, let $P_{j-1}$ be a polynomial
  of degree at most $j-1$. Then there holds
\begin{multline*}
\det_{1\leq i,j\leq n}\left(x_i^{1-j}
\prod_{k=1}^{n-1}(1-c_kx_i)\,P_{j-1}(x_i)
+x_i^{j-1}\prod_{k=1}^{n-1}(1-c_kx_i^{-1})\,
P_{j-1}(x_i^{-1})\right)\\
=2\prod_{i=1}^n P_{i-1}(0)\prod_{1\leq i\le j\leq n-1}(1-c_ic_j)
\prod_{i=1}^n x_i^{1-n}\prod_{1\le i<j\le n}(x_i-x_j)(1-x_ix_j).
\end{multline*}
\end{Corollary}
\begin{proof}
Let $c_n=1$ in Corollary \ref{bdetr1cor} and divide by $\prod_{i=1}^n(1-x_i)$.
\end{proof}

\begin{Remark}\label{wdrem}
If we let $c_j=0$ and $P_j(x)=1$ for all $j$, Corollaries 
\ref{cdetr1cor}, \ref{bdetr1cor} and~\ref{ddetr1cor} reduce, up to
reversing the order of the columns, to the classical Weyl denominator formulas 
\eqref{cwd}, \eqref{bwd} and \eqref{dwd}, respectively. 
Similarly, 
Corollary \ref{apoldet} contains \eqref{awd} as a limit case.
 Thus, Theorems~\ref{adet} and \ref{cdet} give
 elliptic extensions of the
Weyl denominator formulas for the classical root systems. 
\end{Remark}

\section{The Macdonald identities}
\label{secmac}

In Section \ref{secell}, we have focused on the left-hand sides of 
\eqref{mf}, trying to find as general families of $R$ theta functions
as possible, such that the constant $C$ can be determined. We will now
focus on the  right-hand sides,  trying to find a
particularly simple expression for $W_R$ as a determinant. More
precisely,  we want the functions $f_j$ to have  known explicit
Laurent expansions, so that the multiple Laurent 
expansion of $W_R$ can be read off from \eqref{mf}. 

Starting with the case of 
  type $A$, we  observe that the function
 \begin{equation}\label{fm}f_m(x)=x^m\theta((-1)^{n-1}tp^mx^n;p^n),\end{equation}
 with
$m$ an integer, is an $A_{n-1}$ theta function of norm 
$t$. Moreover, its Laurent expansion is known from \eqref{jti}.
Thus, we are led to consider determinants of
 the form
$\det_{ij}\left(f_{m_j}(x_i)\right),$
with $m_j$  integers, hoping that the constant
$$C=\frac{\det_{1\leq i,j\leq
    n}\left(f_{m_j}(x_i)\right)}{\theta(tx_1\dotsm x_n)\,W_{A_{n-1}}(x)} $$
can be evaluated. 

To compute this constant, we specialize the $x_i$ to $n$th
roots of unity, since the theta functions may then be pulled out from
the determinant. To avoid zeroes in the denominator, the $x_i$ should 
be distinct, so we  assume $x_i=\omega^{i-1}$,
with $\omega$  a primitive $n$th root of unity. By the Vandermonde
determinant  \eqref{awd}, we then have
$$\det_{1\leq i,j\leq n}\left(f_{m_j}(\omega^{i-1})\right)
=\prod_{j=1}^n\theta((-1)^{n-1}tp^{m_j};p^n)\prod_{1\leq i<j\leq n}(\omega^{m_j}-\omega^{m_i}). $$   
To obtain a non-trivial result, this should be non-zero, 
so the $m_i$ should be equi\-distributed modulo $n$. Thus, we
 assume $m_i=i-1$. In that case, by \eqref{qrt},
$$\prod_{j=1}^n\theta((-1)^{n-1}tp^{m_j};p^n)=
\theta((-1)^{n-1}t)=
\theta(tx_1\dotsm x_n)\bigg|_{x_i=\omega^{i-1}},$$
which gives
$$\det_{1\leq i,j\leq
  n}\left(x_i^{j-1}\theta((-1)^{n-1}tp^{j-1}x_i^n;p^n) \right)=
\prod_{1\leq i<j\leq
  n}\frac{\omega^{j-1}-\omega^{i-1}}{\omega^{j-1}\theta(\omega^{i-j})}\,W_{A_{n-1}}(x).$$
By \eqref{qr},  the constant simplifies  as  
\begin{multline*} \prod_{1\leq i<j\leq
  n}\frac{\omega^{j-1}-\omega^{i-1}}{\omega^{j-1}\theta(\omega^{i-j})}
= \prod_{1\leq i<j\leq
  n}\frac 1{(p\omega^{j-i},p\omega^{i-j})_\infty}\\
=(p)_\infty^n
\prod_{i,j=1}^n\frac
1{(p\omega^{j-i})_\infty}=(p)_\infty^n\prod_{k=1}^n\frac 1{(p\omega^{k})_\infty^n}
=\frac{(p;p)_\infty^n}{(p^n;p^n)_\infty^n}.
\end{multline*}
Thus,
we arrive at the  $A_{n-1}$ case of Proposition \ref{mdp} below.

For the remaining root systems, we consider the case of Proposition
\ref{wp} when the theta functions are constructed using Lemma
\ref{dl}, with the corresponding functions  $g$  of
the form \eqref{fm}. By similar arguments as for $A_{n-1}$, one is led to
 the following determinants, one for each root system.

\begin{Proposition}\label{mdp}
The following determinant evaluations hold:
$$
\det_{1\leq i,j\leq n}\left(x_i^{j-1}
\theta((-1)^{n-1}p^{j-1}tx_i^n;p^n)\right)=\frac{(p;p)_\infty^{n}}{(p^n;p^n)_\infty^n}
\,\theta(tx_1\dotsm x_n)\,W_{A_{n-1}}(x),
$$
$$
\det_{1\leq i,j\leq n}\left(x_i^{j-n}
\theta(p^{j-1}x_i^{2n-1};p^{2n-1})
-x_i^{n+1-j}
\theta(p^{j-1}x_i^{1-2n};p^{2n-1})
\right)
=\frac{2(p;p)_\infty^{n}}
{(p^{2n-1};p^{2n-1})_\infty^n}\, W_{B_n}(x),$$
$$
\det_{1\leq i,j\leq n}\left(x_i^{j-n-1}
\theta(p^{j-1}x_i^{2n};p^{2n})
-x_i^{n+1-j}
\theta(p^{j-1}x_i^{-2n};p^{2n})
\right)
=\frac{2(p^2;p^2)_\infty(p;p)_\infty^{n-1}}
{(p^{2n};p^{2n})_\infty^n}\,W_{B_n^\vee}(x),
$$
$$
\det_{1\leq i,j\leq n}\left(x_i^{j-n-1}
\theta(-p^{j}x_i^{2n+2};p^{2n+2})
-x_i^{n+1-j}
\theta(-p^{j}x_i^{-2n-2};p^{2n+2})
\right)
=\frac{(p;p)_\infty^{n}}
{(p^{2n+2};p^{2n+2})_\infty^n}\,
W_{C_n}(x),
$$
$$
\det_{1\leq i,j\leq n}\left(x_i^{j-n}
\theta(-p^{j-\frac12}x_i^{2n};p^{2n})
-x_i^{n+1-j}
\theta(-p^{j-\frac12}x_i^{-2n};p^{2n})
\right)
=\frac{(p^{\frac12};p^{\frac12})_\infty(p;p)_\infty^{n-1}}
{(p^{2n};p^{2n})_\infty^n}\,  
W_{C_n^\vee}(x),
$$
$$
\det_{1\leq i,j\leq n}\left(x_i^{j-n}
\theta(-p^{j}x_i^{2n+1};p^{2n+1})
-x_i^{n+1-j}
\theta(-p^{j}x_i^{-2n-1};p^{2n+1})
\right)
=\frac{(p;p)_\infty^n}{(p^{2n+1};p^{2n+1})_\infty^n}\,W_{BC_n}(x)
$$
and, for $n\geq 2$,
$$
\det_{1\leq i,j\leq n}\left(x_i^{j-n}
\theta(-p^{j-1}x_i^{2n-2};p^{2n-2})
+x_i^{n-j}
\theta(-p^{j-1}x_i^{2-2n};p^{2n-2})
\right)
=\frac{4(p;p)_\infty^n}{(p^{2n-2};p^{2n-2})_\infty^n}
\,W_{D_n}(x).
$$
\end{Proposition}

To complete the proof of Proposition \ref{mdp},  all that remains is
to verify the identities for some fixed values of $x_i$. We have
already done this for $A_{n-1}$. In general, we proceed exactly as in 
 \cite{St}. 
Namely, letting $\omega_k$ denote a primitive $k$th root of unity,  we
specialize 
 $x_i$ as 
$x_i=\omega_{4n-2}^{2i-1}$ for $R=B_n$, $x_i=\omega_{4n}^{2i-1}$
for $R=B_n^\vee$,
$x_i=\omega_{2n+2}^{i}$ for $R=C_n$, $x_i=\omega_{2n}^{i}$ for
$R=C_n^\vee$, 
$x_i=\omega_{2n+1}^i$ for $R=BC_n$ and $x_i=\omega_{2n-2}^{i-1}$ for $R=D_n$.
Under these specializations, the theta functions can be pulled out from
the 
determinants, which are then computed by
 the  Weyl denominator formulas \eqref{bwd} (for $B_n$,
$C_n^\vee$ and $BC_n$), \eqref{cwd} (for $B_n^\vee$ and $C_n$) and
\eqref{dwd} (for $D_n$). If we let $Q_R$ denote the quotient of the
determinant and the expression $W_R$, this gives
\begin{align*}
Q_{B_n}&=\frac{\prod_{j=1}^n\theta(-p^{j-1};p^{2n-1})}
{\prod_{j=1}^n(-p\omega_{2n-1}^{\pm(j-n)})_\infty
\prod_{1\leq i<j\leq
  n}(p\omega_{2n-1}^{j-i},p\omega_{2n-1}^{i-j},p\omega_{2n-1}^{i+j-1},p\omega_{2n-1}^{1-i-j})_\infty},\\
Q_{B_n^\vee}&=\frac{\prod_{j=1}^n\theta(-p^{j-1};p^{2n})}
{\prod_{j=1}^n(p^2\omega_{2n}^{\pm(2j-1)};p^2)_\infty
\prod_{1\leq i<j\leq
  n}(p\omega_{2n}^{j-i},p\omega_{2n}^{i-j},p\omega_{2n}^{i+j-1},p\omega_{2n}^{1-i-j})_\infty},\\
Q_{C_n}&=\frac{\prod_{j=1}^n\theta(-p^{j};p^{2n+2})}
{\prod_{j=1}^n(p\omega_{2n+2}^{\pm 2j})_\infty
\prod_{1\leq i<j\leq
  n}(p\omega_{2n+2}^{j-i},p\omega_{2n+2}^{i-j},p\omega_{2n+2}^{i+j},p\omega_{2n+2}^{-i-j})_\infty},\\
Q_{C_n^\vee}&=\frac{\prod_{j=1}^n\theta(-p^{j-\frac12};p^{2n})}
{\prod_{j=1}^n(p^\frac12\omega_{2n}^{\pm j};p^{\frac 12})_\infty
\prod_{1\leq i<j\leq
  n}(p\omega_{2n}^{j-i},p\omega_{2n}^{i-j},p\omega_{2n}^{i+j},p\omega_{2n}^{-i-j})_\infty},\\
Q_{BC_n}&=\frac{\prod_{j=1}^n\theta(-p^{j};p^{2n+1})}
{\prod_{j=1}^n(p\omega_{2n+1}^{\pm j})_\infty(p\omega_{2n+1}^{\pm 2j};p^2)_\infty
\prod_{1\leq i<j\leq
  n}(p\omega_{2n+1}^{j-i},p\omega_{2n+1}^{i-j},p\omega_{2n+1}^{i+j},p\omega_{2n+1}^{-i-j})_\infty},\\
Q_{D_n}&=\frac{2\prod_{j=1}^n\theta(-p^{j-1};p^{2n-2})}
{\prod_{1\leq i<j\leq
  n}(p\omega_{2n-2}^{j-i},p\omega_{2n-2}^{i-j},p\omega_{2n-2}^{i+j-2},p\omega_{2n-2}^{2-i-j})_\infty}.
\end{align*}

It remains  to simplify these expressions
into the form given in Proposition \ref{mdp}.
We indicate a way to organize the computations for $R=B_n$; the other cases
can be treated similarly. We factor
$Q_{B_n}$ as $F_1/F_2F_3$, where 
$$F_1=\prod_{j=1}^n(-p^{j-1};p^{2n-1})_\infty(-p^{2n-j};p^{2n-1})_\infty,
$$
$$F_2=\prod_{j=1}^n(-p\omega_{2n-1}^{j-n})_\infty(-p\omega_{2n-1}^{n-j})_\infty,
 $$ 
$$F_3=\prod_{1\leq i<j\leq
  n}(p\omega_{2n-1}^{j-i},p\omega_{2n-1}^{i-j},p\omega_{2n-1}^{i+j-1},p\omega_{2n-1}^{1-i-j})_\infty.$$

In $F_1$, we
make the change of variables $j\mapsto 2n+1-j$ in the second factor
and use \eqref{qr} to obtain
$$F_1=\prod_{j=1}^{2n}(-p^{j-1};p^{2n-1})_\infty=2(-p;p)_\infty(-p^{2n-1};p^{2n-1})_\infty.$$
Similarly, in $F_2$ we change $j\mapsto 2n-j$ in the second factor, obtaining
$$F_2=\prod_{j=1}^n(-p\omega_{2n-1}^{j-n})_\infty\prod_{j=n}^{2n-1}(-p\omega_{2n-1}^{j-n})_\infty=(-p;p)_\infty(-p^{2n-1};p^{2n-1})_\infty.$$
Finally, in $F_3$ we rewrite the first two factors as
$$\frac 1{(p)_\infty^n}\prod_{i,j=1}^n(p\omega_{2n-1}^{j-i})_\infty. $$
Making the change of variables $i\mapsto 2n-i$, this equals
\begin{equation}\label{ftf}\frac 1{(p)_\infty^n}\prod_{i=n}^{2n-1}\prod_{j=1}^{n}
(p\omega_{2n-1}^{i+j-1})_\infty.\end{equation}
In the fourth factor in $F_3$,
we change $(i,j)\mapsto(n-i,n+1-j)$, which gives
$$
\prod_{1\leq i<j\leq
  n}(p\omega_{2n-1}^{1-i-j})_\infty=\prod_{1\leq j\leq i\leq
  n-1}(p\omega_{2n-1}^{i+j-1})_\infty.$$
Thus, the third and fourth factor can be combined into
$$\prod_{i=1}^{n-1}\prod_{j=1}^{n}
(p\omega_{2n-1}^{j+i-1})_\infty, $$
which, together with \eqref{ftf}, gives
$$F_3=\frac 1{(p)_\infty^n}\prod_{j=1}^{n}\prod_{i=1}^{2n-1}
(p\omega_{2n-1}^{j+i-1})_\infty
=\frac 1{(p)_\infty^n}\prod_{j=1}^{n}
(p^{2n-1};p^{2n-1})_\infty
=\frac{(p^{2n-1};p^{2n-1})_\infty^n}{(p;p)_\infty^n}.
 $$
In conclusion, this shows that
$$Q_{B_n}=\frac{2(p;p)_\infty^{n}}
{(p^{2n-1};p^{2n-1})_\infty^n},$$
in agreement with Proposition \ref{mdp}.

The determinant evaluations in Proposition \ref{mdp} imply 
 the following multiple Laurent expansions. We give two versions of
 each identity, the second one being obtained from the first by an
 application of one of the classical Weyl denominator formulas
\eqref{wd}. To verify that these
 identities agree with Macdonald's, the easiest way is to take the
 second version, replace $p$ by $q$, $m_i$ by $-m_i$ and $x_i$ by $x_i^{-1}$, and
 then compare with how the Macdonald identities are written in \cite{St}. 
(Equation (3.16) in \cite{St} should read
 $c(q)=1/(q)_\infty^n$, not $c(q)=q/(q)_\infty^n$.)

\begin{Corollary}
\label{mlc}
 The following identities hold:
\begin{align*}
(p;p)_\infty^{n-1}\, W_{A_{n-1}}(x)
&=
\sum_{\substack{m_1,\dots,m_n\in\mathbb
  Z\\m_1+\dots+m_n= 0}}\sum_{\sigma\in S_n}\,\sgn(\sigma)
\prod_{i=1}^n
 x_i^{nm_i+\sigma(i)-1}p^{n\binom{m_i}2+(\sigma(i)-1)m_i}\\
&=\sum_{\substack{m_1,\dots,m_n\in\mathbb
  Z\\m_1+\dots+m_n= 0}}
\prod_{i=1}^n
 x_i^{nm_i}p^{n\binom{m_i}2}\prod_{1\leq i<j\leq
  n}(x_jp^{m_j}-x_ip^{m_i})
\end{align*}
\begin{align*}
(p;p)_\infty^n\, W_{B_n}(x)
&=
\sum_{\substack{m_1,\dots,m_n\in\mathbb
  Z\\m_1+\dots+m_n\equiv 0\ (2)}}\sum_{\sigma\in S_n}\,\sgn(\sigma)
\prod_{i=1}^n
 x_i^{(2n-1)m_i}p^{(2n-1)\binom{m_i}2+(n-1)m_i}\\
&\qquad\times
\left( (x_ip^{m_i})^{\sigma(i)-n}-(x_ip^{m_i})^{n+1-\sigma(i)}\right)\\
&=\sum_{\substack{m_1,\dots,m_n\in\mathbb
  Z\\m_1+\dots+m_n\equiv 0\ (2)}}
\prod_{i=1}^n
  x_i^{(2n-1)m_i+1-n}p^{(2n-1)\binom{m_i}2}\\
&\qquad\times\prod_{i=1}^n(1-x_ip^{m_i})\prod_{1\leq i<j\leq
  n}(x_jp^{m_j}-x_ip^{m_i})(1-x_ix_jp^{m_i+m_j}),
\end{align*}
\begin{multline*}\begin{split}
(p^2;p^2)_\infty
(p;p)_\infty^{n-1}\, W_{B_n^\vee}(x)&=\sum_{\substack{m_1,\dots,m_n\in\mathbb
  Z\\m_1+\dots+m_n\equiv 0\ (2)}}\sum_{\sigma\in S_n}\,\sgn(\sigma)
\prod_{i=1}^n
  x_i^{2nm_i}p^{2n\binom{m_i}2+nm_i}\\
&\qquad\times
\left( (x_ip^{m_i})^{\sigma(i)-n-1}-(x_ip^{m_i})^{n+1-\sigma(i)}\right)\end{split}\\
\begin{split}&=\sum_{\substack{m_1,\dots,m_n\in\mathbb
  Z\\m_1+\dots+m_n\equiv 0\ (2)}}
\prod_{i=1}^n
  x_i^{n(2m_i-1)}p^{2n\binom{m_i}2}\\
&\qquad\times
\prod_{i=1}^n(1-x_i^2p^{2m_i})\prod_{1\leq i<j\leq
  n}(x_jp^{m_j}-x_ip^{m_i})(1-x_ix_jp^{m_i+m_j}),
\end{split}\end{multline*}
\begin{align*}
(p;p)_\infty^n\, W_{C_n}(x)
&=
\sum_{m_1,\dots,m_n\in\mathbb Z}\,\sum_{\sigma\in S_n}\,\sgn(\sigma)
\prod_{i=1}^n
  x_i^{(2n+2)m_i}p^{(2n+2)\binom{m_i}2+(n+1)m_i}\\
&\qquad\times
\left( (x_ip^{m_i})^{\sigma(i)-n-1}-(x_ip^{m_i})^{n+1-\sigma(i)}\right)\\
&=\sum_{m_1,\dots,m_n\in\mathbb Z}\,
\prod_{i=1}^n
  x_i^{(2n+2)m_i-n}p^{(2n+2)\binom{m_i}2+m_i}\\
&\qquad\times
\prod_{i=1}^n(1-x_i^2p^{2m_i})\prod_{1\leq i<j\leq
  n}(x_jp^{m_j}-x_ip^{m_i})(1-x_ix_jp^{m_i+m_j})
,
\end{align*}
\begin{multline*}\begin{split}
(p^{\frac12};p^{\frac12})_\infty(p;p)_\infty^{n-1}\, W_{C_n^\vee}(x)
&=
\sum_{m_1,\dots,m_n\in\mathbb Z}\,\sum_{\sigma\in S_n}\,\sgn(\sigma)
\prod_{i=1}^n
  x_i^{2nm_i}p^{2n\binom{m_i}2+(n-\frac12)m_i}\\
&\qquad\times
\left( (x_ip^{m_i})^{\sigma(i)-n}-(x_ip^{m_i})^{n+1-\sigma(i)}\right)
\end{split}\\
\begin{split}&=\sum_{m_1,\dots,m_n\in\mathbb Z}\,
\prod_{i=1}^n
  x_i^{2nm_i+1-n}p^{2n\binom{m_i}2+\frac12m_i}
\\
&\qquad\times\prod_{i=1}^n(1-x_ip^{m_i})\prod_{1\leq i<j\leq
  n}(x_jp^{m_j}-x_ip^{m_i})(1-x_ix_jp^{m_i+m_j}),
\end{split}\end{multline*}
\begin{align*}
(p;p)_\infty^n\, W_{BC_n}(x)
&=
\sum_{m_1,\dots,m_n\in\mathbb Z}\,\sum_{\sigma\in S_n}\,\sgn(\sigma)
\prod_{i=1}^n
  x_i^{(2n+1)m_i}p^{(2n+1)\binom{m_i}2+nm_i}\\
&\qquad\times
\left( (x_ip^{m_i})^{\sigma(i)-n}-(x_ip^{m_i})^{n+1-\sigma(i)}\right)\\
&=\sum_{m_1,\dots,m_n\in\mathbb Z}\,
\prod_{i=1}^n
  x_i^{(2n+1)m_i+1-n}p^{(2n+1)\binom{m_i}2+m_i}
\\
&\qquad\times\prod_{i=1}^n(1-x_ip^{m_i})\prod_{1\leq i<j\leq
  n}(x_jp^{m_j}-x_ip^{m_i})(1-x_ix_jp^{m_i+m_j}),
\end{align*}
\begin{align*}
(p;p)_\infty^n\, W_{D_n}(x)
&=\frac12
\sum_{\substack{m_1,\dots,m_n\in\mathbb
  Z\\m_1+\dots+m_n\equiv 0\ (2)}}\sum_{\sigma\in S_n}\,\sgn(\sigma)
\prod_{i=1}^n
  x_i^{(2n-2)m_i}p^{(2n-2)\binom{m_i}2+(n-1)m_i}\\
&\qquad\times
\left( (x_ip^{m_i})^{\sigma(i)-n}+(x_ip^{m_i})^{n-\sigma(i)}\right)\\
&=\sum_{\substack{m_1,\dots,m_n\in\mathbb
  Z\\m_1+\dots+m_n\equiv 0\ (2)}}
\prod_{i=1}^n
  x_i^{(n-1)(2m_i-1)}p^{(2n-2)\binom{m_i}2}\\
&\qquad\times\prod_{1\leq i<j\leq
  n}(x_jp^{m_j}-x_ip^{m_i})(1-x_ix_jp^{m_i+m_j}),\qquad n\geq 2.
\end{align*}
\end{Corollary}

\begin{proof}
We start  from the determinant evaluations in Proposition \ref{mdp}.  In the cases when there are two theta functions in each
matrix elements (i.e.\ $R\neq A_{n-1}$), we 
apply $\theta(x;p^N)=\theta(p^N/x;p^N)$ to the second one. We then 
expand the left-hand sides using
\eqref{jti}. 
For  $C_n$, $C_n^\vee$ and $BC_n$, this leads  immediately to
 the desired expansions. 

For $A_{n-1}$, expanding also the factor $\theta(tx_1\dotsm x_n)$, we obtain 
\begin{multline*}
\sum_{m_1,\dots,m_n=-\infty}^\infty\sum_{\sigma\in
  S_n}\sgn(\sigma)\prod_{i=1}^n(-1)^{nm_i}
p^{n\binom{m_i}2+(\sigma(i)-1)m_i}t^{m_i}x_i^{nm_i+\sigma(i)-1}\\
=(p)_\infty^{n-1}W_{A_{n-1}}(x)
\sum_{N=-\infty}^\infty(-1)^{N}p^{\binom N2}(tx_1\dotsm x_n)^N.
 \end{multline*}
Viewing this as a Laurent series in $t$, taking the constant term
gives  the desired result. (Picking out any other Laurent coefficient
gives an equivalent identity.)

For   $B_n$,
$B_n^\vee$ and $D_n$, we obtain series with  the right terms but
 different range of summation. More precisely, we find that
$$2X=\sum_{m_1,\dots,m_n\in\mathbb Z}
f(m_1,\dots,m_n), $$
where the identity we wish to prove is
$$X=\sum_{\substack{m_1,\dots,m_n\in\mathbb
  Z\\m_1+\dots+m_n\equiv 0\ (2)}}(-1)^{m_1+\dots+m_n}f(m_1,\dots,m_n)
  $$
 in the cases $B_n$ and $B_n^\vee$,
and 
$$X=\sum_{\substack{m_1,\dots,m_n\in\mathbb
  Z\\m_1+\dots+m_n\equiv 0\ (2)}}f(m_1,\dots,m_n)
  $$
in the case of $D_n$. 
In any case,  it remains to show that
$$\sum_{\substack{m_1,\dots,m_n\in\mathbb
  Z\\m_1+\dots+m_n\equiv 0\ (2)}}f(m_1,\dots,m_n)
=\sum_{\substack{m_1,\dots,m_n\in\mathbb
  Z\\m_1+\dots+m_n\equiv 1\ (2)}}f(m_1,\dots,m_n).$$
To see this, we fix $\sigma$ and  restrict attention to the index
  $m_i$, where
  $i=\sigma^{-1}(1)$. Then, we may write
  $f(m_1,\dots,m_n)=C(g(m_i)+g(m_i+1))$, where $C$ is independent of
  $m_i$ and
\begin{align*}
g(m)&=(-1)^{m}p^{(2n-1)\binom m 2}x_i^{(2n-1)m+1-n},& R&=B_n,\\
g(m)&=(-1)^{m}p^{2n\binom m 2}x_i^{n(2m-1)},& R&=B_n^\vee,\\
g(m)&=p^{(2n-2)\binom m 2}x_i^{(n-1)(2m-1)},& R&=D_n.
\end{align*}
This observation completes the proof.
\end{proof}

\end{document}